\documentclass[11 pt, reqno]{article}

\usepackage{amsmath,amsfonts,amssymb,amsthm}
\usepackage[colorlinks=true,hyperindex=true]{hyperref}
\usepackage{cancel,bbm}
\usepackage{mathabx} 
\usepackage{comment}
\usepackage{enumitem}
\usepackage{cite}

\usepackage{chngcntr}
\counterwithin*{equation}{section}

\usepackage[a4paper, left=2cm, right=2cm, top=2 cm, bottom= 2 cm]{geometry}
\usepackage{color}
\usepackage{fancyhdr}
\usepackage{latexsym}

\newtheorem{ccounter}{ccounter}[section]
\newtheorem{thm}[ccounter]{Theorem}
\newtheorem{lem}[ccounter]{Lemma}
\newtheorem{cor}[ccounter]{Corollary}
\newtheorem{defn}[ccounter]{Definition}
\newtheorem{prop}[ccounter]{Proposition}
\newtheorem{ass}[ccounter]{Assumption}
\newtheorem{ex}[ccounter]{Example}

\def\Xint#1{\mathchoice
   {\XXint\displaystyle\textstyle{#1}}%
   {\XXint\textstyle\scriptstyle{#1}}%
   {\XXint\scriptstyle\scriptscriptstyle{#1}}%
   {\XXint\scriptscriptstyle\scriptscriptstyle{#1}}%
   \!\int}
\def\XXint#1#2#3{{\setbox0=\hbox{$#1{#2#3}{\int}$}
     \vcenter{\hbox{$#2#3$}}\kern-.5\wd0}}

\def\dashint{\Xint-}

\def\bet{\begin{thm}}
\def\eet{\end{thm}}
\def\bel{\begin{lem}}
\def\eel{\end{lem}}
\def\bas{\begin{ass}}
\def\eas{\end{ass}}
\def\bec{\begin{cor}}
\def\eec{\end{cor}}
\def\bed{\begin{defn}}
\def\eed{\end{defn}}
\def\bep{\begin{prop}}
\def\eep{\end{prop}}
\def\beq{\begin{equation}}
\def\eeq{\end{equation}}
\def\proof{\noindent {\bf Proof.}\ \ }
\def\bea{\begin{equation*}}
\def\eea{\end{equation*}}

\def\bex{\begin{ex}}
\def\eex{\end{ex}}

\def\rr{\mathbb{R}}
\def\zz{\mathbb{Z}}

\def\1{\boldsymbol{1}}
\def\Im{\mathrm{Im}}

\def\e{\mathrm{e}}
\def\i{\mathrm{i}}
\def\del{\partial}
\def\d{\mathrm{d}}
\def\eps{\varepsilon}
\renewcommand\leq\varleq
\renewcommand\geq\vargeq
\def\ee{\mathrm{E}}

\def\F{\mathcal{F}}
\def\O{\mathcal{O}}
\def\fa{\mathfrak{a}}
\def\ee{\mathbb{E}}
\def\om{\omega}
\def\pp{\mathbb{P}}

\def\mfa{\mathfrak{a}}

\def\lamN{ \lambda^{(N)}}
\def\lamn{\lambda^{(n)}}
\def\fn{f^{(n)}}
\def\T{\mathcal{T}}
\def\gamK{\gamma^{(K)}}
\def\gamM{\gamma^{(M)}}
\def\muN{\mu^{(N)}}

\def\A{\mathcal{A}}
\def\rhosc{\rho_\mathrm{sc}}
\def\msc{m_{\mathrm{sc}}}
\def\E{\mathcal{E}}

\def\Fxi{\mathcal{F}_\xi }
\def\hati{\hat{i}}
\def\bx{{\boldsymbol{x}}}
\def\by{{\boldsymbol{y}}}
\def\R{\mathcal{R}}
\def\G{\mathcal{G}}
\def\Gxi{\mathcal{G}_\xi}
\def\hatx{\hat{x}}
\def\haty{\hat{y}}
\def\tilx{\tilde{x}}
\def\tily{\tilde{y}}
\def\gamc{\gamma^{(c)}}
\def\gamN{\gamma^{(N)}}
\def\gamM{\gamma^{(M)}}
\def\gamx{\gamma^{(x)}}
\def\gamy{\gamma^{(y)}}
\def\nux{\nu^{(x)}}
\def\nuy{\nu^{(y)}}
\def\S{\mathcal{S}}
\def\US{\mathcal{U}^{(\mathcal{S})}}
\def\UA{\mathcal{U}^{(\mathcal{A})}}
\def\mfm{\mathfrak{m}}

\begin{document}
\title{Deformed GOE}


\begin{table}
\centering

\begin{tabular}{c}
\multicolumn{1}{c}{\Large{\bf Edge scaling limit of Dyson Brownian motion}}\\
\\
\multicolumn{1}{c}{\Large{\bf at equilibrium with general $\beta \geq 1$}}\\
\\
\\
\end{tabular}
\begin{tabular}{c }
Benjamin Landon\\
\\
 Department of Mathematics \\
 \small{MIT} \\
\small{blandon@mit.edu}  \\
\\
\end{tabular}
\\
\begin{tabular}{c}
\multicolumn{1}{c}{\today}\\
\\
\end{tabular}

\begin{tabular}{p{15 cm}}
\small{{\bf Abstract:} For general $\beta \geq 1$, we consider Dyson Brownian motion at equilibrium and prove convergence of the extremal particles to an ensemble of continuous sample paths in the limit $N \to \infty$.  For each fixed time, this ensemble is distributed as the Airy$_\beta$ random point field.  We prove that the increments of the limiting process are locally Brownian.  When $\beta >1$ we prove that after subtracting a Brownian motion, the sample paths are almost surely locally $r$-H{\"o}lder for any $r<1-(1+\beta)^{-1}$.  Furthermore for all $\beta \geq 1$ we show that the limiting process solves an SDE in a weak sense.  When $\beta=2$ this limiting process is the Airy line ensemble. }
\end{tabular}
\end{table}

\section{Introduction and main results}

Given an infinite sequence of independent two-sided standard Brownian motions $\{ B_i \}_{i=1}^\infty$ we consider for each integer $N \geq 1$, the solution $ \{ \lamN_i (t) \}_{i=1}^N$ to the following system of stochastic differential equations,$^1${\let\thefootnote\relax\footnotetext{1. We refer the interested reader to, e.g., \cite{AGZ} for a pedagogical discussion of the existence and uniqueness of strong solutions of Dyson Brownian motion.}} for $\beta \geq 1$,
\begin{align} \label{eqn:DBM-def}
\d \lamN_i (t) = \sqrt{ \frac{2}{ \beta} } \d B_i (t) + \sum_{ j \neq i }  \frac{1}{ \lamN_i (t) - \lamN_j (t) } \d t - \frac{ \lamN_i (t)}{2 N^{1/3} } \d t, \quad t > -N
\end{align}
where the initial data $\{ \lamN_i (-N) \}_{i=1}^N$ is distributed as the Gaussian $\beta$-ensemble scaled so that it is invariant under the above dynamics (so that the extremal particles are near $\pm2N^{2/3}$).  The above system is called Dyson Brownian motion (DBM), and the invariant Gaussian $\beta$-ensemble  has the density on $\rr^N$ given by,
\beq \label{eqn:gbe}
p^\beta (\mu_1, \dots, \mu_N) \d \mu  := \frac{1}{Z} \exp \left[ - \beta N^{-1/3} \sum_{i=1}^N \frac{\mu_i^2}{4} + \beta \sum_{1 \leq i < j \leq N } \log (\mu_j - \mu_i) \right] \1_{\{ \mu_1 < \mu_2 < \dots < \mu_N\}} \d \mu.
\eeq
For example, under this scaling, the quantity $(-2N^{2/3} - \lambda_1 (t))$ converges to the  Tracy-Widom$_\beta$ distribution as $N \to \infty$.  

Note that the $B_i (t)$'s do not depend on $N$, so that the systems $\{ \lamN_i (t) \}_i$ for different $N$ have the same Brownian motions.  Using this coupling we prove the following estimate.
\bet \label{thm:mr}
There is an $\fa >0$ and a $N_0 >0$ so that the following holds.  Let $N$ and $M$ satisfy $M \geq N \geq N_0$.  Then there is an event $\F_{N, M}$ such that $\pp [ \F_{N, M} ] \geq 1 - N^{-\fa}$ on which,
\beq
\sup_{1 \leq i \leq N^{\fa}} \sup_{ |t| \leq N^{\fa}}  \left| (\lambda^{(N)}_i (t) + 2N^{2/3} )  - ( \lambda^{(M)}_i (t) + 2 M^{2/3} ) \right| \leq N^{-\fa}.
\eeq
\eet
This has some obvious consequences.  Define $N_k = 2^k$.  Clearly, the following limit exists almost surely,
\beq \label{eqn:lam-inf-def}
\lambda_i (t) := \lim_{k \to \infty}  (\lambda_i^{(N_k)} (t) + 2(N_k)^{2/3} ).
\eeq
Furthermore, we see that almost surely, $\lambda_i^{(N_k)} (t)$ converges uniformly on compact subsets $[-T, T]$ to $\lambda_i (t)$ so the latter are continuous and finite functions of $t$.  Additionally, the $\lambda_i(t)$ stay ordered (of course, they may not be \emph{strictly} ordered in the limit $N \to \infty$).  

 The estimate of Theorem \ref{thm:mr} immediately implies that the processes $\lambda_i^N(t)$ converge in probability uniformly on compact sets to the limiting process $\lambda_i (t)$.  We summarize these observations in the following.
\bet \label{thm:conv-prob}
Let $\{ \lambda_i (t) \}_{i \geq 1, t \in \rr}$ be defined as in \eqref{eqn:lam-inf-def}.  Then each sample path is continuous in $t$ and for all $i \geq 1$ we have,
\beq
\lambda_i (t) \leq \lambda_{i+1} (t)
\eeq
for all $ t\in \rr$.  Furthermore, there is a $\fa >0$ and $ C_0 >0$ so that,
\beq
\pp \left[ \sup_{1 \leq i \leq N^{\fa} ,|t| \leq N^{\fa}} | ( \lamN_i (t) + 2 N^{2/3} )- \lambda_i (t) | \geq N^{-\fa} \right] \leq C_0 N^{-\fa}.
\eeq
\eet
We will use the polynomial rate of convergence to deduce some properties of the limiting system from the finite dimensional ones.  Using the fact that the finite systems satisfy \eqref{eqn:DBM-def} as well as estimates for the Gaussian $\beta$-ensemble we prove the following theorem.
\bet \label{thm:local-brown}
For any fixed $i$, $s$ and $S$, the quantity
\beq \label{eqn:brown-inc}
\frac{1}{ \sqrt{ \eps}} \sup_{  0 \leq t \leq S} \left| ( \lambda_i (s+ \eps t) - \lambda_i (s) ) - (2 / \beta)^{1/2} (B_i (s + \eps t) - B_i (s) ) \right|
\eeq
converges to $0$ in probability as $\eps \to 0$.  If $\beta > 1$, then for any $r$ satisfying,
\beq
r < 1 - \frac{1}{ 1+ \beta}
\eeq
there are exponents $\mfm_1 >0$ and $\mfm_2 >0$ and a constant $C_1 >0$ so that for any $\eps >0$, we have for $i \leq \eps^{-\mfm_1}$ that the estimate,
\beq \label{eqn:brown-glob}
\pp\left[ \sup \left\{ \frac{ \left| \lambda_i (t) - \lambda_i (s) - (2/\beta)^{1/2} (B_i (t)-B_i (s) ) \right| }{ |t-s|^r}: |t|, |s| \leq 1 , |t-s| \leq \eps \right\} > C_1  \right] \leq C_1 \eps^{\mfm_2}.
\eeq
holds.  In particular, the path $t \to \lambda_i (t) - (2/\beta)^{1/2} B_i (t)$ is almost surely locally $r$-H{\"o}lder for any $r < 1 -(\beta+1)^{-1}$.
\eet

A natural question about the limiting process is  whether it solves any sort of limiting version of the SDE \eqref{eqn:DBM-def}.  Towards this, we prove the following. 
\bet  \label{thm:sde}
There is a $C>0$ so that the following holds.  For integer $K \geq 1$ and for $i \leq K^{1/{100}}$ 
\begin{align}
\pp \bigg[ \sup_{0 \leq t \leq K^{1/100} }  \bigg| \lambda_i (t) - \lambda_i (0) -& (2/\beta)^{1/2} (B_i (t) - B_i (0) ) \notag\\
&- \sum_{ j \neq i, j \leq K}  \int_0^t \frac{1}{ \lambda_i (s) - \lambda_j (s) } \d s + a K^{1/3} t \bigg| > K^{-b} \bigg] \leq C K^{-b}
\end{align}
where $b = 1/100$ and $a$ is the constant $a= (16/(3 \pi^2 ))^{1/3}$. 
\eet

\subsection{Relation to other work}

Dyson Brownian motion \eqref{eqn:DBM-def} appears as the flow of eigenvalues of a Brownian motion on the space of symmetric (resp., Hermitian, quaternionic self-dual) in the case $\beta=1$ (resp., $\beta=2, 4$) \cite{dyson}, and the measure \eqref{eqn:gbe} is the joint distribution of the eigenvalues of the invariant Gaussian Orthogonal (resp., Unitary, Symplectic) Ensemble.

There has been much work studying the existence of scaling limits of solutions to \eqref{eqn:DBM-def} in both the bulk and the edge, as well as various  characterizations of the limiting process.  In the case $\beta=2$ the correlation functions of \eqref{eqn:DBM-def} have a special determinantal structure.  For this determinantal case, convergence of the correlation functions of \eqref{eqn:DBM-def} in both the bulk and edge for general cases of initial data was obtained by Katori-Nagao-Tanemura \cite{KNT}.  Moreover, in the equilibrium $\beta=2$ case, the limiting object is known as the Airy line ensemble \cite{CH}, which is believed to be the universal scaling limit of various two-dimensional statistical mechanical models; we will discuss further this aspect momentarily.

The works of Osada \cite{O1,O2}, Kawamoto-Osada \cite{KO} and Osada-Tanemura \cite{OT} take a different point of view of the edge limit of the system \eqref{eqn:DBM-def}.  In the cases of $\beta=1, 2, 4$ they have shown existence and uniqueness of an infinite dimensional system of SDEs for general classes of initial data related to \eqref{eqn:DBM-def} (its definition requires care as there is cancelation between the interaction term and the final confining term of \eqref{eqn:DBM-def} which are both formally infinite in the $N \to \infty$ limit).  Moreover, they have also shown convergence of the finite dimensional solutions to the solution of the infinite system.  
Note that while we handle general $\beta \geq 1$, we treat only the stationary case and do not consider any general class of initial data.  Moreover, our only result on any sort of limiting system of SDEs is the weak form of Theorem \ref{thm:sde}.

The work of Tsai \cite{T} considered the bulk limit of \eqref{eqn:DBM-def} and proved the convergence to a system of infinite-dimensional SDEs for general $\beta \geq 1$ (note that the bulk has a different scaling than that indicated by our set-up of \eqref{eqn:DBM-def}).  This result holds both for equilibrium initial data and a class of non-equilibrium initial data.

In the work \cite{Sodin}, Sodin considered the eigenvalues of submatrices of time-varying Wigner matrices.  The eigenvalues of a matrix and those of its top left $k \times k$ minors comprise an interlacing triangular array, each level being the eigenvalues of a minor.  Sodin proved that under some conditions on the matrix process, that these time-varying triangular arrays have  a limit at the spectral edge.   In the example of Dyson Brownian motion, the top level of his limiting process agrees with ours in the case of $\beta=1, 2, 4$ (his work explicitly only deals with $\beta=1, 2$ but the modification to $\beta=4$ is likely straightforward).  In the $\beta=2$ case this process was previously studied by Adler-Nordenstam-Van Moerbeke \cite{ANvM} and Ferrari-Frings \cite{FF}.

As can be seen from the discussion, prior to our work, dynamical convergence of the edge scaling limit of Dyson Brownian motion at equilibrium was unavailable outside of the classical $\beta=1, 2, 4$.  For general $\beta$ only the convergence of the fixed time distributions was known; this was proven by Ramirez, Rider and Vir{\'a}g \cite{RRV}.  They moreover offered an attractive description of the limiting distribution as the smallest eigenvalues of the Stochastic Airy operator, a random Schr{\"o}dinger operator on the half-line, confirming the prediction of Edelman-Sutton \cite{ES}. 
  It is a motivating question on whether the dynamical limit may be described as the eigenvalues of a time-dependent random Schr{\"o}dinger operator.  In this direction, matrix models for Dyson Brownian motion have been studied by Holcott-Paquette \cite{HP} and Allez-Bouchard-Guionnet \cite{ABG,AG}.  An alternative approach to the edge limits of $\beta$-ensembles was developed by Gorin-Shkolnikov in \cite{gorin2018stochastic}. 

In a different direction, the works of Gorin and Shkolnikov have introduced and studied multilevel systems of Dyson Brownian motion  \cite{GS1,GS2}.  That is, they consider for general $\beta$, time dependent triangular arrays whose fixed-time distributions are a generalization (in the parameter $\beta$) of the GUE corners process (i.e., the eigenvalues of a GUE matrix and its the top left $k\times k$ corners), and each level of which is a Dyson Brownian motion.  For $\beta \geq 4$ they prove the joint convergence of the (gaps between) the left-most particle of each level of the DBM (whereas we consider the left-most particles of a single level of DBM).

Very recently, Gorin and Kleptsyn \cite{gorin2020universal} studied the convergence of the edge limit of Dyson Brownian motion in the limit $\beta=\infty$ (note that there is a  rescaling of the particles by $\sqrt{\beta}$ before taking the limit $\beta \to \infty$).  Their result \cite[Theorem 1.2]{gorin2020universal} can be viewed as an integrable analogue of our Theorem \ref{thm:conv-prob}.  It would be interesting to take the $\beta \to \infty$ limit of our process and try to recover their limit; that is, to prove that the $\beta, N \to \infty$ limits can be interchanged.

In the case $\beta=2$ the limiting process is known as the Airy line ensemble \cite{CH}.  This ensemble is believed to be a universal scaling limit for various two-dimensional statistical mechanics models.  For example,  Pr{\"a}hofer and Spohn \cite{PS} showed that it arises as the scaling limit a polynuclear growth model.  We refer the reader to the works \cite{DV,DMV, Q, C-KPZ} and the references therein for more on the scaling limits of statistical mechanical models.

With regards to the Airy line ensemble, a recent research direction has been the investigation of the local Brownian nature of the sample paths.  The work of Corwin and Hammond proved a Brownian Gibbs property for the Airy line ensemble and consequently obtained that the sample paths are locally absolutely continuous with respect to Brownian motion.  This point was quantified in terms of estimates on the Radon-Nikodyn derivative by Calvert-Hammond-Hegde \cite{CHH}.  We were able to use the explicit rate of convergence provided by Theorem \ref{thm:conv-prob} together with the fact that the finite-dimensional system satisfies the SDEs of \eqref{eqn:DBM-def} to obtain the results of Theorem  \ref{thm:local-brown} on the local properties of the sample paths of our limit.  

Our method of proving Theorem \ref{thm:conv-prob} is partly based on prior work of the author together with Sosoe and Yau \cite{edgeDBM,landonyau,fixed} on the local ergodicity of Dyson Brownian motion with general initial data.  The motivation for that work is the study of universality of eigenvalues of random matrices; we refer to the works  \cite{erdos2017dynamical,localspectral,tao2010random,tao2011random,erdos2019matrix} and the references therein.  The application of Dyson Brownian motion in the context of universality for general $\beta$ originates with the works of Erd{\"o}s-Schlein-Yau \cite{localrelaxation}. 

The work \cite{edgeDBM} establishes the local ergodicity of DBM with general initial data at the spectral edge in the case of $\beta=1, 2, 4$.  The result proven there is that given two solutions to \eqref{eqn:DBM-def} with different initial data at $t=0$, the difference between the two solutions (up to a deterministic shift and scaling) after any polynomially large time $t=N^{\delta}$ is $o(1)$ for particles near the edge.  An adaptation of this method was used by Cipolloni-Erd{\H{o}}s-Kr{\"u}ger-Schr{\"o}der in their work on the cusp universality of random matrices \cite{cusp1,cusp2}.  For $\beta=2$ the limiting point process is known as the Pearcey process and has a determinantal structure; however for $\beta=1$ no such process exists.  In \cite{cusp2}, they applied their adaptation of the work \cite{edgeDBM} to compare random matrix systems of different sizes, allowing for the proof of existence of a limiting $\beta=1$ point process.  We use a similar idea here and couple all of the systems of different sizes together.

The analysis of \cite{edgeDBM} requires as input certain rigidity estimates, which are high probability estimates on the particle locations at a near-optimal scale.  Rigidity estimates are usually established using random matrix methods involving analysis of the resolvent.   This is the primary reason why the results of \cite{edgeDBM} are resticted to $\beta=1, 2, 4$ where random matrix methods are available.  We remark that the optimal rigidity results in the cusp case were proven in \cite{cusp1}, being a key technical estimate in their proof of universality. 

By working at equilibrium we bypass the random matrix restriction $\beta=1, 2, 4$ as rigidity results for general $\beta$-ensembles were proven by Bourgade-Erd{\H{o}}s-Yau \cite{bourgade2012bulk,bourgade2014universality,bourgade2014edge}.  However, the method of \cite{edgeDBM} requires the use of certain interpolating systems between the two sets of initial data under consideration.  If we were to apply this method to the set-up considered here, the interpolating systems would not be in equilibrium and we would have no rigidity results.

Rigidity for the process \eqref{eqn:DBM-def} for non-equilibrium classes of initial data was established in the bulk by the author together with Huang \cite{HL} and at the edge by Adhikari-Huang \cite{AH} (we note that these results apply to a more general situation than that considered here, where the quadratic potential implicit in \eqref{eqn:DBM-def} is replaced by a general confining potential).  However, the edge rigidity is delicate, and the interpolating processes constructed in \cite{edgeDBM} do not satisfy the assumptions of \cite{AH}.  Instead of attempting to modify the work \cite{AH}, we give a different approach to analyzing the system \eqref{eqn:DBM-def} than \cite{edgeDBM}; this approach avoids the use of any interpolating ensembles at the cost of somewhat weaker estimates.  This method may be independent interest, applicable in situations where rigidity can be established via other means but the use of interpolating ensembles needs to be avoided.  A drawback of this approach is that one is required to establish level repulsion estimates, which estimate the probability that two particles get too close.  In our situation, these estimates at the edge were basically established by Bourgade-Erd{\H{o}}s-Yau in \cite{bourgade2014edge} (building on ideas of \cite{Gap,bourgade2014universality}).  Additionally, similar to \cite{edgeDBM}, our approach uses the important coupling idea of \cite{homogenization} as well as the energy estimate of \cite{bourgade2014edge}.  

We would also like to mention that Bourgade has recently given a new approach to the ergodicity of Dyson Brownian motion \cite{bourgade2018extreme}. As an application he obtained the first explicit rate of convergence to the Tracy-Widom distribution for general Wigner matrices.   

As a final side comment, we point out that using the universality results of  \cite{AH, edgeDBM}, it is possible to deduce that our limiting process is also the edge scaling limit of (up to deterministic re-scalings):
\begin{enumerate}
\item The Langevin dynamics of $\beta$-ensembles with general potential at equilibrium.
\item The ``long-time'' limit of Langevin dynamics of $\beta$-ensembles with general potential and certain classes of deterministic initial data.  To be more precise,  if one considers the solution $\muN_i(t)$ of \eqref{eqn:DBM-def} with initial data at $t=0$ obeying the conditions of \cite{AH}, then for any $\omega >0$, the limit of $\{ \mu_i (t + N^{\omega} \}_{i, t}$ will coincide with our limiting process $\{ \lambda_i (t) \}_{i, t}$. 
\end{enumerate}
In light of the above, the rough heuristic is then that if a solution of \eqref{eqn:DBM-def} (or its general potential analogue) with some set of initial data has the condition that the fixed time distributions converge to the Airy$_\beta$ random point field, then  it is expected that the multi-time distribution converges to that of $\{ \lambda_i (t) \}_{i, t}$.


\vspace{5 pt}

\noindent{\bf Acknowledgements.}  The author thanks Vadim Gorin for suggesting this problem as well as comments on a draft of this work.  The author thanks Amol Aggarwal for comments on a draft and for suggesting further lines of investigation of the limiting process, as well as Philippe Sosoe for enlightening discussions.  The author also thanks Paul Bourgade for guidance concerning  level repulsion estimates.

\subsection{Organization of remainder of paper}

In Section \ref{sec:gbe} we collect some known results about the Gaussian $\beta$-ensemble.  In Section \ref{sec:conv} we prove the results regarding convergence of DBM edge,  Theorems \ref{thm:mr} and \ref{thm:conv-prob}.  In Section \ref{sec:local} we prove our results on the properties of the limiting process, Theorems \ref{thm:local-brown} and \ref{thm:sde}.  An appendix collects a technical level repulsion result which is proven by following the  argument in \cite{bourgade2014edge}.

\section{Auxilliary results for the Gaussian $\beta$-ensemble} \label{sec:gbe}
 
We will denote expectation and probability with respect to the Gaussian $\beta$-ensemble $p^\beta (\mu_1, \dots, \mu_N)$ defined in \eqref{eqn:gbe} by $\ee^\beta$ and $\pp^\beta$.  As is well-known, the empirical measure of the particles under $p^\beta$ when rescaled by $N^{2/3}$ converges to Wigner's semicircle distribution,
\beq
\rhosc (E) \d E := \1_{ \{ |E \leq 2 \} } \frac{1}{ 2 \pi} \sqrt{4 - E^2} \d E.
\eeq
Denote by $\gamN_i$ the $N$-quantiles of the semicircle distribution,
\beq
\frac{i}{N} = \int_{-2}^{\gamN_i } \rhosc (E) \d E.
\eeq
At one point we will we have use of the Stieltjes transform,
\beq
\msc (z) = \int \frac{ \rhosc (x) }{ (x-z) } \d z = \frac{-z + \sqrt{z^2-4}}{ 2}.
\eeq
We let $\Fxi$ be the event,
\beq
\Fxi = \bigcap_{i=1}^N \left\{ | \mu_i - N^{2/3} \gamN_i | \leq N^{\xi} (\hati )^{-1/3} \right\},
\eeq
where
\beq
\hati := \min \{ i, N+1-i \}.
\eeq
The following rigidity result was proven in Theorem 2.4 of \cite{bourgade2014edge}.  
\bet \label{thm:rig}
For any $\xi >0$ there are constants $c>0$ and $C>0$ so that,
\beq
\pp^\beta[ \Fxi ] \geq 1 - C \e^{-N^{c}}.
\eeq
\eet
We require a level repulsion estimate which controls the probability that neighbouring particles get too close under the Gaussian $\beta$-ensemble measure.  For the most part, the result we need was proven in Theorem 3.2 of \cite{bourgade2014edge}.  However, the result stated there is not quite the full statement we require, as it misses a few particles near the edge.  Fortunately, the proof technique applies with almost no modifications to give the result we need; we give the details in Appendix \ref{a:lr} where the proof of the following is provided, essentially copying the work done in \cite{bourgade2014edge}. 
\bet \label{thm:lr}
For all $\eps >0$ and $r>0$ there is a constant $C>0$ we have for all sufficiently small $\xi$ (depending on $\eps, r$) that for all $s>0$,
\beq
\pp^\beta [ \Fxi \cap \{  \mu_{i+1} - \mu_i \leq s \hati^{-1/3} ]  \leq C N^{\eps} s^{1+\beta-r},
\eeq
for all $i \leq N^{9/10}$.
\eet

\section{Convergence calculations} \label{sec:conv}

Given an infinite sequence of two-sided Brownian motions $\{ B_i \}_{i=1}^\infty$ we consider for each $N$, the solution $\lamN_i (t)$ to the SDEs, for $\beta \geq 1$,
\begin{align}
\d \lamN_i (t) = \sqrt{ \frac{2}{ \beta} } \d B_i (t) + \sum_{ j \neq i }  \frac{1}{ \lamN_i (t) - \lamN_j (t) } \d t - \frac{ \lamN_i (t)}{2 N^{1/3} } \d t, \quad t > -1
\end{align}
where the initial data $\{ \lamN_i (-1) \}_{i=1}^N$ is distributed as a Gaussian $\beta$-ensemble scaled so that the largest particle is at $2N^{2/3}$.  This distribution is stationary for the above system.  Note that the $B_i (t)$'s do not depend on $N$, so that the systems $\{ \lamN_i (t) \}_i$ for different $N$ have the same Brownian motions.

Fix $N$ sufficiently large and $M$ satisfying $N \leq M \leq 10 N$.  We define,
\beq
x_i (t) = \lamN_i (t) + 2N^{2/3}, \qquad y_i (t) = \lambda^{(M)}_i (t) +2M^{2/3}.
\eeq
They satisfy the following equations,
\begin{align}
\d x_i (t) = \sqrt{ \frac{2}{ \beta}} \d B_i (t) + \sum_{j \neq i } \frac{1}{ x_i(t) -x_j (t) } \d t + N^{1/3} \d t - \frac{x_i (t)}{2 N^{1/3}} \d t,
\end{align}
and
\begin{align}
\d y_i (t) = \sqrt{ \frac{2}{ \beta}} \d B_i (t) + \sum_{j \neq i } \frac{1}{ y_i (t) - y_j (t) } \d t + M^{1/3} \d t - \frac{y_i (t) }{ 2 M^{1/3} } \d t.
\end{align}
The goal of this section is to prove the following theorem
\bet \label{thm:te}
There is a $\mfa >0$ so that for all sufficiently large $N$ the following holds with probability at least $1-N^{-\mfa}$,
\beq
\sup_{1 \leq i \leq N^{\mfa}} \sup_{ |t| \leq N^{\mfa}} \left| x_i (t)   - y_i (t)  \right| \leq N^{-\mfa}
\eeq
\eet
Theorem \ref{thm:mr} is a straightforward consequence of this and a dyadic argument.

\noindent{\bf Proof of Theorem \ref{thm:mr}}.   Given $N$ and $M$ define a sequence $\{ n_k\}_{k=0}^m$ s.t. $n_0 = N$, $n_m = M$ and for $0 < i < m$, 
\beq
n_i = 2^i 2^{ \lfloor \log_2 (N) \rfloor }
\eeq
with $m-1 = \lfloor \log_2 (M) \rfloor - \lfloor \log_2(N) \rfloor$.  Then we write,
\beq
(\lambda^{(M)}_i (t) + 2 M^{2/3} )  - ( \lamN_i (t) +2N^{2/3} )= \sum_{j=1}^m  \left[ \lambda^{(n_j)}_i (t) + 2 n_j^{2/3} ) - ( \lambda_i^{(n_{j=1} ) }  (t) + 2 n_{j-1}^{2/3} ) \right].
\eeq
Let $\E_j$ be the event,
\beq
\sup_{ |t| \leq N^{\fa} }  \sup_{ 1 \leq i  \leq N^{\fa} }  \left|  \lambda^{(n_j)}_i (t) + 2 n_j^{2/3} ) - ( \lambda_i^{(n_{j=1} ) }  (t) + 2 n_{j-1}^{2/3} )  \right| \leq (n_{j-1})^{-\fa}.
\eeq
Since $N \leq n_j$ for all $j \geq 1$, by Theorem \ref{thm:te} we have that $\pp[ \E_j ] \geq 1 - (n_{j-1} )^{-\fa}$.  We let $\F_{N, M} = \bigcap_{j=1}^m \E_j$.  The claim follows from the fact that
\beq
\sum^m_{j=1} (n_j)^{-\fa} \leq C 2^{ - \fa \lfloor \log_2(N) \rfloor } \sum_{j>0} 2^{-j} \leq C N^{-\fa}.
\eeq
\qed

\subsection{Quantiles}
Let $\gamN$ and $\gamM$ be the $N$- and $M$-quantiles of $\rhosc$,
\beq
\frac{i}{N} = \int_{-2}^{\gamN_i} \d \rhosc (x) \d x, \qquad \frac{j}{M} = \int_{-2}^{\gamM_i} \d \rhosc (y ) \d y.
\eeq
By direct calculation for $i \leq \alpha N$ for small $\alpha$,
\beq
N^{2/3} ( \gamN_i + 2 ) = \left( \frac{2}{3} i \pi \right)^{2/3} + \O \left( \frac{i^{5/3}}{N} \right)
\eeq
Define now the quantiles,
\beq
\gamx_i := N^{2/3} ( \gamN_i + 2), \qquad \gamy_i := M^{2/3} ( \gamM_i + 2),
\eeq
and the measures,
\beq
\nux (E) \d E := N^{1/3} \rhosc(N^{-2/3} E-2 ) \d E, \qquad \nuy (E) \d E = M^{1/3} \rhosc (M^{-2/3} E-2 ) \d E. 
\eeq

For $\xi >0$ let us denote the event, 
\begin{align}
\G_\xi &:= \{ |x_i (t) - \gamx_i | \leq N^{\xi} N^{-2/3} \hati^{-1/3} : 1 \leq i \leq N, -N \leq t \leq N \} \nonumber\\
 &\cap   \{ |y_i (t) - \gamy_i| \leq M^{\xi} M^{-2/3} \hati^{-1/3} : 1 \leq i \leq N, -N \leq t \leq N \}
\end{align}
We have the following lemma, a similar form of which appeared before in \cite{fixed}.   It may be deduced from \cite{AH}.  We provide a self-contained proof in Appendix \ref{a:sc}. 
\bel \label{lem:stoch-cont}
For any $\xi$, $D>0$ we have that
\beq
\pp [ \G_\xi ] \geq 1 - N^{-D}
\eeq
for large enough $N$.
\eel

\subsection{Regularized dynamics}

We first introduce the following regularization, similar to \cite{homogenization}.  Fix,
\beq
C_r = 10^6,
\eeq
and define
\beq
\eps_{jk} = \begin{cases} N^{-C_r}, &  j > k  \\ - N^{-C_r}, & j < k \end{cases}.
\eeq
We also take a $t_0$ satisfying,
\beq
-N \leq t_0 \leq 0.
\eeq
We allow $t_0$ to depend on $N$.
We define the regularizations $\hatx_i$ by
\beq
\d \hatx_j (t) = \sqrt{ \frac{2}{ \beta}} \d B_j (t) + \sum_{k \neq j } \frac{1}{ x_j (t) - x_k (t) + \eps_{jk} }\d t +N^{1/3} \d t - \frac{x_j (t)}{2 N^{1/3}} \d t , \qquad t > t_0
\eeq
and
\beq
\hatx_j (t_0) = x_j (t_0).
\eeq
We construct the regularized dynamics $\haty_j$ similarly, by
\beq
\d \haty_j (t) = \sqrt{ \frac{2}{ \beta}} \d B_j (t) + \sum_{k \neq j } \frac{1}{ y_j (t) - y_k (t) + \eps_{jk} }\d t +M^{1/3} \d t - \frac{y_j (t)}{2 M^{1/3}} \d t , \qquad t > t_0
\eeq
and
\beq
\haty_j (t_0) = y_j (t_0).
\eeq
The following provides an estimate for the effect of the regularization.  The proof is similar to \cite{homogenization}.
\bel \label{lem:reg}
There is an event $\F_1$ with $\pp [ \F_1 ] \geq 1 - N^{-100}$ on which,
\beq
| \hatx_j (t) - x_j (t) |+ | \haty_j (t) - y_j (t) |  \leq N^{-10}, \qquad t_0 \leq t \leq N, 1 \leq j \leq N^{1/2}.
\eeq
\eel
\proof The difference satisfies,
\beq
\d ( \hatx_j (t) - x_j (t) )= \sum_{k \neq j } \frac{ \eps_{jk} }{ (x_j - x_k + \eps_{jk} )(x_j -x_k) } \d t.
\eeq
We have, for sufficiently small $\xi >0$,
\begin{align}
&\sum_{ j \leq N^{1/10}, k \neq j } \ee\left[ \1_{\Gxi}  \int_{-N}^N \left|   \frac{ \eps_{jk} }{ (x_j (t)- x_k(t) + \eps_{jk} )(x_j (t)-x_k(t)) } \right| \right] \nonumber \\
\leq & N^{1-C_r /2 } \sum_{j \leq N^{1/10} } \int_{-N}^N \ee\left| \1_{\Gxi} \frac{1}{ |x_j(t) - x_{j+1} (t) |^{3/2} } \right| \nonumber \\
\leq & C N^{-10^5}.
\end{align}
In the last line we used Theorem \ref{thm:lr}.  The estimate for the difference $\hatx_j - x_j$ follows from integration and Markov's inequality.  The estimate for $\haty_j - y_j$ follows analogously. \qed

\subsection{Approximate cut-off dynamics}
We introduce two exponents, 
\beq
\om_K>0, \qquad \tau_1 >0
\eeq
and the parameters
\beq
K = N^{\om_K}, \qquad t_1 = N^{\tau_1}.
\eeq
We will assume,
\beq
\om_K < 10^{-1}, \qquad \tau_1 < \om_K / 10^3.
\eeq
We require a further scale given by $\delta_c >0$.  We assume,
\beq
\delta_c < 10^{-3},
\eeq
and define,
\beq
\gamc := \gamx_{\lfloor K+K^{\delta_c}\rfloor }.
\eeq
We now define some process $\{ \tilx_j (t)\}_{ 1 \leq j \leq K}$ and $\{ \tily_j (t) \}_{1 \leq j \leq K}$ for $t \geq t_0$, which we will refer to as the approximate cut-off dynamics.  We define the $\tilx_i's$ as the solution to,
\beq \label{eqn:xtil-def}
\d \tilx_j (t) = \sqrt{ \frac{2}{ \beta}} \d B_j (t) + \sum_{ k \neq j,  k \leq K} \frac{1}{ x_j (t) - x_k (t) + \eps_{jk}} \d t + \1_{ \{ x_j (t) \leq \gamc\} } \left( \int_{ \gamc}^\infty \frac{ \nux (x))}{ x_j (t) - x } \d x \right) \d t + N^{1/3} \d t , 
\eeq
for $t > t_0$ and 
with initial data $\tilx_j (t_0 ) = \hatx_j (t_0) = x_j (t_0)$.  Similarly, we define
\beq \label{eqn:ytil-def}
\d \tily_j (t) = \sqrt{ \frac{2}{ \beta}} \d B_j (t) + \sum_{ k \neq j,  k \leq K} \frac{1}{ y_j (t) - y_k (t) + \eps_{jk}} \d t +\1_{ \{ y_j (t) \leq \gamc\} }  \left(\int_{\gamc}^\infty \frac{ \nux (x) }{ y_j (t) - x } \d x \right) \d t + N^{1/3} \d t , 
\eeq
for $t >t_0$ and 
with initial data $\tily_j (t_0) = \haty_j (t_0) = y_j (t_0)$.   A few remarks are in order.  First, in the dynamics for $\tily$ we use some of the $x$ quantities, i.e.,  the measure $\nux$ and $N^{1/3}$ (instead of $M^{1/3}$ and $\nuy$).  We will provide momentarily a few deterministic calculations which will be used to account for this difference when we estimate the effect of our cut-offs.  

Secondly, note  that on the event $\Gxi$ for $\xi >0$ sufficiently small, the indicator function in the above dynamics is identically $1$.

We now complete the aforementioned deterministic calculations. 
\bel
Suppose $a$ satisfies,
\beq
- \gamx_{2K} \leq a \leq \gamx_{K + K^{\delta_c} - K^{\delta_c/4}}.
\eeq
Then,
\begin{align}
\int_{\gamc}^\infty \frac{ \nux (x) }{ a - x } \d x = \dashint_{0}^{\gamc} \frac{ \sqrt{x}}{\pi (a-x) } \d x -N^{1/3}  + \sqrt{ (a)_- }+ \O( N^{-1/3} \log(N) |\gamc|^{3/2} ),
\end{align}
where the integral on the RHS is interpreted as a principal value.
\eel
\proof 
We write,
\begin{align}
\int_{\gamc}^\infty \frac{ \nux (x) }{ a - x } \d x &= \dashint_{-\infty}^\infty \frac{ \nux (x) }{ a - x } \d x -  \dashint_{\gamc}^\infty \frac{ \nu(x) }{ a - x } \d x.
\end{align}
The first term equals,
\begin{align}
\dashint_{-\infty}^\infty \frac{ \nux (x) }{ a - x } \d x &= - N^{1/3} \msc ( N^{-2/3} a - 2 ) = -N^{1/3} + \sqrt{ (a)_- } + \O\left( |a|  N^{-1/3}  +|a|^{3/2} N^{-2/3} \right).
\end{align}
For the second term, we start with
\begin{align}
\nux (x) = N^{1/3} \rhosc (N^{-2/3} x - 2 ) =   \frac{1}{\pi} \sqrt{x} - \frac{1}{2 \pi} N^{-2/3} x^{3/2} \frac{1}{ \sqrt{4} + \sqrt{4 - N^{-2/3} x } }.
\end{align}
Define $g(x) := (2 + \sqrt{4-N^{-2/3} x } )^{-1}$.
If $a \leq 0$ then,
\beq
\int_{0}^{\gamc}  N^{-2/3} \left| \frac{x^{3/2} g(x)}{ x-a} \right| \d x \leq C N^{-2/3} |\gamc|^{3/2}.
\eeq
Fix some $\eta >0$, with $ \eta < \gamc / 2$.   If $0< a < \eta$,
\beq \label{eqn:det-a1}
\dashint_0^{\gamc}  \frac{ x^{3/2} g(x) }{ x-a} \d x = \dashint_{0}^{2a} \frac{x^{3/2} g(x)}{ x -a } \d x + \int_{2a}^{\gamc} \frac{x^{3/2} g(x) }{ x -a } \d x.
\eeq
The first term satisfies,
\beq
 \dashint_{0}^{2a} \frac{x^{3/2} g(x)}{ x -a } \d x =  \int_{0}^{2a} \frac{x^{3/2} g(x) - a^{3/2} g(a)}{ x -a } \d x  = \O \left( |a|^{3/2} \right)
\eeq
The second term satisfies,
\beq
\int_{2a}^{\gamc}  \left|\frac{x^{3/2} g(x) }{ x -a }  \right| \d x \leq C \int_{2a}^{\gamc} |x|^{1/2} \d x \leq C |\gamc|^{3/2}.
\eeq
If $a > \eta$, we break the integral up into three segments.
\begin{align}
\dashint_0^{\gamc}  \frac{ x^{3/2} g(x) }{ x-a} \d x  = \int_{0}^{a - \eta} +\dashint_{a-\eta}^{a+\eta} + \int_{a+\eta}^{\gamc} \frac{ x^{3/2} g(x) }{ x-a} \d x
\end{align}
The first one is estimated by,
\beq \label{eqn:det-a2}
\left| \int_{0}^{a- \eta}  \frac{ x^{3/2} g(x) }{ x-a} \d x \right| \leq |a|^{3/2} C \int_{0}^{a-\eta} \frac{1}{ |x-a| } \d x \leq C |a|^{3/2} |\log (\eta ) |.
\eeq 
  The second is estimated similarly to \eqref{eqn:det-a1} and is $\O( \eta |a|^{1/2} )$.  The third is estimated similarly to \eqref{eqn:det-a2} and is $\O ( |\gamc|^{3/2} | \log ( \eta ) )$. We take $\eta = N^{-20}$. \qed

\bel \label{lem:xy-rig-dif}
Suppose $a$ satisfies,
\beq
- \gamx_{K} \leq a \leq \gamx_{K + K^{\delta_c} - K^{\delta_c/2}}.
\eeq
Then,
\begin{align}
\left| \int_{\gamx_{K + K^{\delta_c}}}^\infty \frac{ \nux (x) }{a-x} \d x +N^{1/3} - M^{1/3} - \int_{\gamy_{K+K^{\delta_c}}}^\infty \frac{ \nuy(x)}{ a -x } \d x \right| \leq C N^{-1/3} \log(N) K .
\end{align}
\eel
\proof
Note that since $\gamx_j = \gamy_j + \O\left(N^{-1} K^{5/3} \right) $ for $j \leq 2K$, the condition on $a$ implies that,
\beq
-\gamy_{2K} \leq a \leq \gamy_{K+K^{\delta_c} - K^{\delta_c/4} }.
\eeq  
Hence we can apply the previous lemma with the $y$-quantities in place of $x$ (i.e., $\nuy$ and $M^{1/3}$ instead of $\nux$ and $N^{1/3}$).  We see that in order to conclude the proof  we need to estimate,
\beq
\left| \int_{\gamx_{K+K^{\delta_c}}}^{ \gamy_{K+K^{\delta_c}}}\frac{\sqrt{x}}{x - a} \right|.
\eeq
The condition on $a$ implies that the denominator is greater than $1$.  Hence, the integral is bounded by,
\beq
\left| \int_{\gamx_{K+K^{\delta_c}}}^{ \gamy_{K+K^{\delta_c}}}\frac{\sqrt{x}}{x - a} \right| \leq C|\gamc|^{1/2} K^{5/3} N^{-1} 
\eeq
and we conclude the proof from our assumption that $K\leq N^{1/10}$ and from the estimate $|\gamc| \leq C K^{2/3}$.
\qed

We often have to control quantities like $i^{2/3}  - j^{2/3}$, so  the following elementary lemma is useful.
\bel
For $a$ and $b$ positive,
\beq
|a^{2/3} - b^{2/3} | \asymp \frac{|a-b|}{a^{1/3}+b^{1/3}}.
\eeq
\eel
\proof This follows simply from the identity,
\beq
A^3 - B^3 = (A-B)(A^2 +AB + B^2) = (A^2 - B^2) \frac{A^2 + AB + B^2}{A+B}
\eeq
and taking $A = a^{1/3}$ and $B= b^{1/3}$. \qed

We will use the above lemma without comment throughout this section.

We can now prove the following estimate controlling the approximate cut-off dynamics introduced above.
\bel \label{lem:cutoff-1}
Let $\delta >0$.  There is a $c>0$ so that the following holds.  There is an event $\F_2$ of probability greater than $\pp [\F_2] \geq 1 - K^{-c \delta}$, on which it holds for every $1 \leq j \leq K$ that,
\beq
\sup_{ 0 \leq t \leq t_1} | \tilx_j (t_0 + t) - \hatx_j (t_0 + t)  |  + | \tily_j (t_0+t) - \haty_j (t_0 + t) | \leq C (1+t_1) \frac{K^{1/3} K^{\delta+\delta_c}}{K-j+1}
\eeq
\eel
\proof We have,
\begin{align}
\d (\hatx_j - \tilx_j  ) (t)  &=\left(  \sum_{K < l \leq K + K^{\delta_c} } \frac{1}{ x_j (t) - x_l (t) + \eps_{jl} } \right) \d t \\
&+ \left( \sum_{l \geq K + K^{\delta_c}} \frac{1}{ x_j (t) - x_l (t) + \eps_{jl} } - \1_{ \{ x_j (t) \leq \gamc \} } \int_{\gamc }^\infty  \frac{  \d \nux (x) }{ x_j (t) - x } \d x \right) \d t\\
&+ \frac{1}{2 N^{1/3}}( x_j  ) (t) ) \d t \\
& =: A_{1, j} (t) \d t+A_{2, j} (t)  \d t+ A_{3, j} (t) \d t
\end{align}
If $j \leq K - K^{\delta/10}$, then we see that on $\G_\xi$ for $\xi$ sufficiently small, that 
\beq
\left| A_{1, j} (t)\right|  \leq C \frac{ K^{\delta_c}}{ (K+1)^{2/3} - j^{2/3}} \leq C \frac{ K^{\delta_c} K^{1/3}}{K - j +1}
\eeq
for all $t \in [-1, 1]$.  For $K - K^{\delta/10} \leq j \leq K$ we estimate,
\beq
\int_{0}^{t_1} | A_{j, 1}(t_0 +s ) | \d s \leq \sum_{K - K^{\delta/10} \leq k \leq K} \int_{0}^{t_1} |A_{k, 1} (t_0 + s) |\d s.
\eeq
We have,
\beq
 \ee\left[ \1_{\Gxi}    \sum_{K - K^{\delta/10} \leq k \leq K} \int_{0}^{t_1} |A_{k, 1} (t_0 + s) |\d s \right] \leq C K^{\delta_c+2 \delta} K^{1/3} t_1
\eeq
Hence by Markov's inequality, there is a single event of probability at least $1- K^{-c \delta}$ on which for every $1 \leq j \leq K$ we have,
\beq
\int_0^{t_1} |A_{1, j} (t_0 +s ) | \d s \leq (1 + t_1) \frac{ K^{3 \delta + \delta_c} K^{1/3}}{K - j + 1}.
\eeq
On $\Gxi$ for sufficiently small $\xi$ we have, that 
\begin{align}
|A_{2, j} | &\leq \sum_{ l \geq K + K^{\delta_c}} \left|  \frac{1}{ x_j (t) - x_l (t) + \eps_{j, l} } - \int_{\gamx_l}^{\gamx_{l+1}} \frac{ \d \nux (x)}{ x_j (t) - x } \right|\nonumber\\ 
& \leq C \sum_{ l \geq K + K^{\delta_c} }  \frac{N^{\xi}}{ (j^{2/3} - l^{2/3} )^2 \min \{ l^{1/3}, (N+1 - l)^{1/3} \} } \nonumber\\
& \leq C N^{\xi} \sum_{l \geq K + K^{\delta_c}} \frac{l^{1/3}}{ (j-l)^2} + C N^{-1} \leq  C  \frac{K^{1/3+\delta}}{K-j+1}.
\end{align} 
We assume that $\xi$ is small enough that $N^{\xi } \leq K^{\delta}$. 
For the last inequality we used,
\begin{align}
\sum_{ l > K+1} \frac{l^{1/3}}{( j-l)^2} \leq C \sum_{l >K+1} \frac{ j^{1/3}}{ (j-l)^2} + \frac{1}{ (j-l)^{5/3}} \leq C \frac{K^{1/3}}{ (K+1-j) } .
\end{align}
On the event $\Gxi$ we for sufficiently small $\xi >0$ that
\beq
|A_{3, j} | \leq C \frac{K}{N^{1/3}}.
\eeq
This completes the estimate for the difference $\hatx_j - \tilx_j$.  For the difference $\haty_j - \tily_j$ there is an additional term,
\beq
A_{4, j} :=  \int_{\gamx_{K + K^{\delta_c}}}^\infty \frac{ \nux (x) }{a-x} \d x +N^{1/3} - M^{1/3} - \int_{\gamy_{K+K^{\delta_c}}}^\infty \frac{ \nuy(x)}{ a -x } \d x.
\eeq
This is handled by Lemma \ref{lem:xy-rig-dif}. \qed

For the gaps we have a better estimate.
\bel \label{lem:cutoff-2}
Let $\delta >0$.  There is an event $\F_3$ with probability at least $\pp [\F_3] \geq 1- K^{-c \delta}$ on which,
\beq
\sup_{ 0 \leq t_1 \leq t } \left| ( \tilx_a - \tilx_b ) - ( \hatx_a - \hatx_b ) \right| (t_0 + t) \leq K^{\delta+\delta_c} (1+t_1)  \frac{|a-b|}{a^{1/3} + b^{1/3}} \frac{K^{2/3}}{ (K-a+1)(K-b+1) }
\eeq
for all $a, b \leq K$.  The same estimate holds for the quantities $(  \tily_a - \tily_b ) - ( \haty_a - \haty_b )$.  Moreover the same estimate holds with the $\hatx_i$'s replaced by $x_i$'s, etc.
\eel
\proof We write,
\begin{align}
\d \left( ( \tilx_a - \tilx_b ) - ( \hatx_a - \hatx_b ) \right) &= \left( \sum_{ K \leq l \leq K + K^{\delta_c}} \frac{ x_a - x_b }{ (x_l - x_a + \eps )( x_l - x_b + \eps ) } \right) \\
+ (x_a - x_b) & \left( \sum_{l \geq K + K^{\delta_c} } \frac{1}{(x_a - x_l - \eps ) (x_b - x_l - \eps )} - \int_{\gamc}^\infty \frac{ \d \nux (x)}{ (x_a -x )(x_b -x ) } \right) \d t \\
&+ \frac{ x_a - x_b}{2 N^{1/3}} \d t =: \left( D_1 (t) + D_2 (t) + D_3 (t) \right) \d t.
\end{align}
As in the previous lemma,$|D_3 (t) | \leq C K N^{-1/3}$ on $\Gxi$.  We next estimate $D_1$.  We have,
\beq
|D_1| \leq K^{\delta_c} \frac{ |x_a-x_b|}{ (x_{K+1} - x_a + \eps )(x_{K+1} -x_b + \eps ) }
\eeq
On $\Gxi$ for $\xi$ sufficiently small, we have for $p \leq K$,
\beq
(x_{K+1}-x_p + \eps ) \geq ((K+1)^{2/3} - p^{2/3} )K^{1/3} K^{-\delta} (x_{K+1} - x_K + \eps ).
\eeq
Hence, for all $a, b \leq K$ we have on $\Gxi$,
\beq
\int_{0}^{t_1} |D_1 (t_0+s) | \d s \leq K^{4 \delta + \delta_c} \frac{|a-b|}{a^{1/3} + b^{1/3}} \frac{K^{2/3}}{ (K-a+1)(K-b+1) } \int_{t_0}^{t_0+t_1} \int \frac{1}{ K^{2/3} (x_{K+1} (s)- x_{K} (s)+ \eps )^2} \d s.
\eeq
By Theorem \ref{thm:lr} and Markov's inequality, there is an event of probability at least $1- K^{-c \delta}$ on which,
\beq
\int_{t_0}^{t_0+t_1} \frac{1}{ K^{2/3} (x_{K+1} (s)- x_{K} (s)+ \eps )^2} \leq (1+t_1) K^{\delta}.
\eeq
For $D_2$ we have by rigidity that on $\Gxi$ for small enough $\xi >0$, for $a \leq b$ we have,
\begin{align}
\frac{1}{ x_b -x_a} \left|D_2\right| \leq K^{\delta} \sum_{l \geq K + K^{\delta_c}} \frac{l^{-1/3}}{ (l^{2/3} - a^{2/3} )(l^{2/3} -b^{2/3} )^2 } \leq C K^{\delta} \sum_{l \geq K + K^{\delta_c}} \frac{l^{2/3}}{ (l-a) (l-b)^2 }.
\end{align}
We estimate the final sum by,
\begin{align}
\sum_{l \geq K + K^{\delta_c}} \frac{l^{2/3}}{ (l-a) (l-b)^2 } & \leq C b^{2/3} \sum_{l \geq K + K^{\delta_c}} \frac{1}{ (l-a)(l-b)^2 }+ C \sum_{l \geq K +  K^{\delta_c}} \frac{1}{ (l-a)(l-b)^{4/3}} \nonumber\\
&\leq C \frac{ b^{2/3}}{ (K-a+1) (K-b+1) } +C \frac{1}{(K-a+1)(K-b+1)^{1/3} } \nonumber\\
&\leq C \frac{K^{2/3}}{ (K-a+1)(K-b+1)}.
\end{align}
This completes the estimate for the differences of the $\hatx_i$'s and $\tilx_i$'s.  The estimates for the $\haty_i$'s and $\tily_i$'s are similar. \qed 

\subsection{Parabolic equation for difference of cut-off dynamics}

We now consider the differences,
\beq
u_j (t) := \tilx_j (t) - \tily_j (t), \qquad 1 \leq j \leq K.
\eeq
On the event $\Gxi$ for sufficiently small $\xi >0$, this function obeys the following discrete parabolic equation
\begin{align}
\frac{\d }{\d t} u_i (t) =  \sum_{j \neq i, j \leq K} B_{ij} (t) (u_j (t) - u_i (t) ) - W_i (t) u_i (t) + \xi_i (t)+ q_i (t), 
\end{align}
with coefficients,
\begin{align}
B_{ij} (t) &:= \frac{1}{ (x_i (t) - x_j (t) + \eps_{ij} )(y_i (t) - y_j (t) + \eps_{ij} ) }, \nonumber\\
 W_i (t) &:= \1_{ \{ x_i (t) \leq \gamc\} } \1_{ \{ y_i (t) \leq \gamc \} } \int_{\gamc}^\infty \frac{ \d \nux (x) }{ (x_i (t) - x)(y_i (t) - x)} \d x,
\end{align}
and forcing terms,
\begin{align}
\xi_i (t) &:= \sum_{j \neq i, j \leq K} B_{ij} (t) ( [(\tilx_i - \tilx_j ) - ( x_i -x_j ) ] - [ (\tily_i - \tily_j ) - (y_i -y_j ) ] ) \nonumber\\
q_i(t) &:= W_i (t)( \tilx_i - x_i  - (\tily_i - y_i ) ),
\end{align}
where we omitted the argument $(t)$ from some of the $\tilx_i$'s, $x_i$'s, etc. 
Note that we restrict to the event $\Gxi$ so that the terms $W_i (t)$ have a simple form; recall that the indicator functions appearing in \eqref{eqn:xtil-def} and \eqref{eqn:ytil-def} are identically $1$ on $\Gxi$.  We introduce the operator $\A$ by
\beq \label{eqn:A-def}
( \A (t) w)_i := \sum_{ j \neq i, j \leq K } B_{ij}(t) (w_j - w_i ) - W_i (t) w_i .
\eeq

\subsection{Estimates of forcing terms}
We now estimate the quantities $\xi_i (t)$ and $q_i(t)$ which appeared in the parabolic equation that we derived in the previous subsection. 
\bel \label{lem:xi}
For any $\delta >0$ there is an event $\F_4$ with probability $\pp[ \F_4 ] \geq 1 - K^{-c \delta}$ on which we have for all $1 \leq a \leq K$ and all $0 \leq t \leq t_1$,
\beq
 | \xi_a (t_0 + t) | \leq (1 + t_1) \frac{ |B_{a, a+1}(t_0 + t) | + \1_{\{ a>1 \} } |B_{a, a-1} (t_0+t) | }{a^{2/3}} \frac{K^{1+\delta+\delta_c}}{(K-a+1)^2}.
\eeq
\eel
\proof On the event $\F_3 \cap \Gxi$ where $\F_3$ is the event described in Lemma \ref{lem:cutoff-2} and $\xi>0$ small enough we have,
\begin{align}
|\xi_a| & \leq (1+t_1) \sum_{|j-a| \leq K^{\delta}} \frac{ \left( |B_{a, a+1} | +  \1_{\{ a>1 \} } |B_{a, a-1} | \right) K^{\delta+\delta_c} K^{2/3}}{ a^{1/3} (K-a+1)^2} \nonumber\\
&+ (1+t_1) \sum_{ |j-a| > K^{\delta}} \frac{ a^{1/3} + j^{1/3}}{ |j-a|} \frac{K^{2/3} K^{\delta+\delta_c}}{ (K-a+1)(K-j+1) }  \nonumber\\
&\leq (1+t_1) \frac{|B_{a, a+1} | +  \1_{\{ a>1 \} } |B_{a, a-1} | }{a^{2/3}} K^{3 \delta+\delta_c} \frac{K}{(K-a+1)^2} \nonumber\\
&+ (1+t_1) \sum_{ |j-a| > K^{\delta}} \frac{ a^{1/3} + j^{1/3}}{ |j-a|} \frac{K^{2/3} K^{\delta+\delta_c}}{ (K-a+1)(K-j+1) }.
\end{align}
Note that in the above, all summations are restricted to indices $j \leq K$. The same remark holds for the following calculation in which we
 estimate the sum in the last line by,
\begin{align}
& \sum_{|j -a| > K^{\delta}} \frac{a^{1/3} + j^{1/3}}{ |a-j| } \frac{K^{2/3}}{ (K-a +1)(K-j + 1)} \nonumber \\
\leq&  \frac{K}{K-a+1} \sum_{ j \leq (K+a)/2, j \neq a} \frac{1}{ |a-j| (K-j+1) }+ \frac{K}{K-a+1} \sum_{ j > (K+a)/2} \frac{1}{|a-j| (K-j+1)} \nonumber \\
\leq & \frac{2K}{K-1+a} \sum_{j \leq (K+a)/2, j \neq a} \frac{1}{ |a-j| (K-a+1) } + \frac{2K}{K-a+1} \sum_{j > (K+a)/2} \frac{1}{ (K-a+1)(K-j+1) } \nonumber \\
\leq & C \frac{K \log(K) }{ (K-a+1)^2}
\end{align}
Hence, on the event $\F_3 \cap \Gxi$ we see that,
\beq
|\xi_a| \leq (1+t_1) \frac{|B_{a, a+1} | +  \1_{\{ a>1 \} } |B_{a, a-1} | }{a^{2/3}}  \frac{K^{1+5 \delta+\delta_c}}{(K-a+1)^2}.
\eeq
This is the desired estimate. \qed

\bel \label{lem:qi}
Let $\delta >0$.  There is an event $\F_5$ with probability at least $\pp [ \F_5] \geq 1 - K^{-c \delta}$ on which for every $i$,
\beq
\sup_{0 \leq t \leq t_1} |q_i (t_0 + t) | \leq (1+t_1) \frac{ K^{1+\delta+\delta_c}}{(K-i+1)^2}.
\eeq
\eel
\proof On the event $\F_2 \cap \Gxi$, where $\F_2$ is the event of Lemma \ref{lem:cutoff-1} we have for sufficiently small $\xi$,
\begin{align}
|q_i| &\leq (1+t_1) \frac{K^{1/3} K^{\delta+\delta_c}}{K-i+1} \sum_{j > K } \frac{j^{2/3}}{ (j-i)^2} \leq (1+t_1) \frac{K^{1/3+\delta+\delta_c}}{K-i+1} \sum_{j > K } \frac{i^{2/3}}{ (j-i)^2} + \frac{1}{(j-i)^{4/3}} \nonumber\\
&\leq  (1+t_1) \frac{K^{1/3+\delta+\delta_c}}{K-i+1} \left( \frac{K^{2/3}}{ K-i+1} + \frac{1}{ (K-i+1 )^{1/3}} \right) \leq  C(1+t_1) \frac{K^{1+\delta+\delta_c}}{(K-i+1)^2}.
\end{align}
This is the claim. \qed

\subsection{Finite speed of propogation}

We now split our operator into a short-range and long-range part.  We fix a cut-off $\ell$,
\beq
\ell := K^{2/3+\eps_\ell}
\eeq
where
\beq
0 < \eps_\ell < 1/6.
\eeq
We then define,
\beq
\A = \S + \R
\eeq
where 
\beq
(\S u)_j :=  \sum_{ |j-k | \leq \ell, j \neq k} B_{jk} (u_k - u_j )  - W_j u_j,
\eeq
and
\beq
( \R u)_j :=  \sum_{ k : |k-j | > \ell } B_{jk} (u_k - u_j ).
\eeq
On $\Gxi$ for $\xi>0$ sufficiently small, we have
\beq
\sum_{k : |j-k | > \ell} \left| B_{jk}(t) \right| \leq \frac{C}{K^{\eps_\ell}}.
\eeq
Therefore, for every $1 \leq p \leq \infty$,
\beq
||\R||_{\ell^p \to \ell^p } \leq \frac{C}{K^{\eps_\ell}}.
\eeq
We denote the semigroups of the operators $\A$ and $\S$ by $\UA(s, t)$ and $\US(s, t)$, respectively so that (for example) for $w_0 \in \rr^K$, the function
\beq
 w(t) := \UA(s, t) w_0
\eeq
solves 
\beq
\frac{ \d }{ \d t} w(t) = \A (t) w (t), \quad t > s, \qquad w(s) = w_0.
\eeq
It is standard that both $\UA(s, t)$ and $\US (s, t)$ are contractions on any $\ell^p$ space which implies,
\beq
\sum_j \left| \US_{ij} (s, t) \right| \leq 1, \qquad \sum_i \left| \US_{ij} (s, t) \right| \leq 1,
\eeq
and similarly for $\UA_{ij} (s, t)$.   Moreover, the matrix elements of $\US$ and $\UA$ are all positive, 
\beq
0 \leq \UA_{ij} (s, t), \qquad 0 \leq \US_{ij} (s, t)
\eeq
 The following exponential decay estimate is a modification of a similar proof in \cite{Gap}.   The range in \eqref{eqn:expdecay-range} is not optimal but suffices for our purposes.  
\bel \label{lem:expdecay}
Let $\delta >0$.  There is a $c>0$ so that the following holds.  There is an event $\F_6$ with probability at least $\pp[\F_6] \geq 1- K^{-c \delta}$ on which, the following estimate holds for all $a, b$ satisfying
\beq \label{eqn:expdecay-range}
|a - b | \geq K^{\delta} K^{5/6} \sqrt{1+t_1}
\eeq
and all $t_0 \leq s \leq t \leq t_0  +t_1$:
\beq
|\US_{ab} (s, t) | \leq C \e^{ - K^{c \delta}}.
\eeq
\eel
\proof Let $\theta \geq \ell$.  Define $r(t) = \US (s, t) \delta_b$, and $f(t)$ by,
\beq
f(t) = \sum_{j=1}^K \phi_j r_j (t)^2, \qquad \phi_j = \exp \left[ \theta^{-1} |j-b| \right].
\eeq
We have, following \cite{Gap},
\begin{align}
f'(t) &= 2 \sum_j \phi_j \sum_{k : |j-k| \leq \ell} r_j (t) B_{kj} (t) (r_k - r_j ) (t) - 2 \sum_j \phi_j r_j^2 W_j \nonumber\\
&\leq 2 \sum_j \phi_j \sum_{k : |j-k| \leq \ell} r_j (t) B_{kj} (t) (r_k - r_j ) (t)  \nonumber \\
&= \sum_{ |j - k | \leq \ell } B_{kj} (t) (r_k -r_j )(t) [ r_j (t) \phi_j - r_k (t) \phi_k ] \nonumber\\
&= \sum_{ |j-k | \leq \ell } B_{kj} (r_k -r_j ) (t) \phi_j (r_j - r_k ) (t) \nonumber\\
&+ \sum_{ |j - k | \leq \ell} B_{kj} (t) (r_k -r_j ) (t) [ \phi_j - \phi_k ] r_k (t). \label{eqn:fe-1}
\end{align}
For the second term, we have by the Schwarz inequality,
\begin{align}
\sum_{ |j - k | \leq \ell} B_{kj} (t) (r_k -r_j ) (t) [ \phi_j - \phi_k ] r_k (t) &\leq \frac{1}{2} \sum_{ |j-k | \leq \ell } B_{kj} (r_k - r_j )^2 \phi_k \nonumber\\
&+ \frac{1}{2} \sum_{ |j-k| \leq \ell} B_{kj} (t) r_k^2 \phi_k^{-1} (\phi_k - \phi_j )^2.
\end{align}
The term on first line on the RHS is absorbed into the term in the second-last line of \eqref{eqn:fe-1}, this latter term being negative.
Assuming that $\ell \leq \theta$,
\beq
\phi_k^{-2} [ \phi_j - \phi_k ]^2\leq C \frac{ |j-k|^2}{ \theta^2} , \qquad |j-k | \leq \ell.
\eeq
Therefore,
\begin{align}
f'(t) &\leq C \sum_{ |j-k| \leq \ell} B_{kj} \phi_k^{-1} ( \phi_j  - \phi_k)^2 r_k^2 \nonumber\\
& \leq \frac{C}{\theta^2} B_{kj}|j-k|^2 \phi_k r_k (t)^2 \nonumber\\
& \leq C \theta^{-2} \sum_{K^{\delta} \leq |j-k| \leq \ell} B_{kj} |j-k|^2 \phi_k r_k^2 +C \theta^{-2}\sum_{|j-k| \leq K^{\delta}} B_{kj} |j-k|^2 \phi_k r_k^2 \nonumber \\
&= A_1 + A_2.
\end{align}
For the term $A_1$ we have on $\Gxi$ for sufficiently small $\xi >0$ that,
\beq
A_1 \leq  C \theta^{-2} K^{2/3} \ell \sum_k \phi_k r_k^2 \leq \theta^{-2} K^{5/3} f(t).
\eeq
The other term is bounded by
\beq
A_2 \leq \theta^{-2} \left( \sum_{ j, k : |j-k| \leq K^{2\delta}} B_{kj} K^{2 \delta} \right) f(t).
\eeq
By Gronwall, we therefore have that for all $t_0 \leq 0 \leq s \leq t \leq t_0 + t_1$,
\beq
f(t) \leq \exp \left[ C \theta^{-2} \left( (1+t_1) K^{5/3} + \int_{0}^{t_1} \sum_{|j-k| \leq K^{\delta}} B_{kj}(t_0+s) K^{2 \delta} \right) \right] f(s).
\eeq
Note that $f(s) =1$ by definition.  
Now we have,
\beq
\ee\left[ \int_{0}^{t_1} \sum_{|j-k| \leq K^{\delta} } B_{kj} K^{2 \delta} \right] \leq (1+t_1) K^{4\delta} \sum_{j=1}^K \sup_{ 0 \leq t \leq t_1} \ee |B_{j, j+1}(t_0+t) | \leq (1+t_1) K^{5 \delta}  K^{5/3}.
\eeq
By Markov's inequality there is an event of probability at least $1- K^{-\delta}$ on which,
\beq
\int_0^{t_1}  \sum_{|j-k| \leq K^{\delta} } B_{kj} K^{2 \delta}  \leq (1+t_1) K^{6 \delta} K^{5/3}. 
\eeq
Note that this event does not depend on the specific choice of $b$ or $s$. 
Hence, on the intersection of this event and $\Gxi$ we have that $f(t) \leq C$ (for any choice of $b$ or $s$) as long as
\beq
\theta = K^{10\delta} K^{5/6} \sqrt{1+t_1}.
\eeq
Note our assumptions on $\ell$ imply $\ell \leq K^{5/6}$ so $\ell \leq \theta$ is satisfied with the above choice. 
The claim follows. \qed

Fix a large natural number $m$.  We have,
\begin{align}
\UA (s, t) &= \US (s, t) \\
&+ \sum_{k=1}^{m-1} \int_{s \leq s_1 \leq \dots \leq s_k t} \US(s_k, t) \R (s_k ) \US (s_{k-1}, s_k ) \dots \R (s_1 ) \US (s, s_1 ) \d s_1 \dots \d s_k \nonumber \\
&+ \int_{s \leq s_1 \leq \dots \leq s_m \leq t } \UA (s_m, t) \R (s_m) \dots \R (s_) \US (s, s_1 ) \d s_1 \dots \d s_{m} \nonumber \\
&=: \US(s, t) +\sum_{k=1}^{m-1} \int_{s \leq s_1 \leq \dots s_k \leq t } A^{(k)} \d s_1 \dots \d s_k + \int_{s \leq s_1 \leq \dots \leq s_m \leq t } B^{(m)} \d s_1 \dots \d s_m.
\end{align}
The operator $\R$ appears $m$ times in $B^{(m)}$ and so for any $1 \leq p \leq \infty$,
\beq
||B^{(m)} ||_{\ell^p \to \ell^p } \leq \frac{C}{ K^{m \eps_\ell} }.
\eeq
We now bound $A^{(k)}$.
\bel
Let $m$ be as above.  Let $\delta >0$.   There is an event with probability at least $\pp[ \F_7] \geq 1- K^{- c \delta}$ on which the following holds.  For any $p \geq K^{ \delta} K^{5/6} \sqrt{1+t_1} $ we have for $k \leq m$,
\beq
| A^{(k)}_{ap} | \leq \frac{C_m}{ p^{4/3}}, \qquad a  \leq p \frac{m-k +1}{200 m}  =: Q_k
\eeq
for some constant $C_m$.
\eel
\proof We will repeatedly use the estimates,
\beq \label{eqn:USest}
\sum_{j} |\US_{ij} (s, t) | \leq 1, \qquad \sum_{i} | \US_{ij} (s, t) | \leq 1.
\eeq
We define,
\beq
D_\S := K^{\delta/2} K^{5/6} \sqrt{1 + t_1}.
\eeq
We let $\F_7$ be the intersection of event of Lemma \ref{lem:expdecay} that exponential decay holds for $|a-b| \geq D_\S$ and $\Gxi$ for $\xi$ sufficiently small.

The proof of the lemma is by induction on $k$.  For $k=1$ we have
\beq
A^{(1)}_{ap } = \sum_{i, j} \US_{ai} \R_{ij} \US_{jp}
\eeq
By the exponential decay estimates we can restrict the summation to $|i-a| \leq D_\S$ and $|j-p| \leq D_\S$.  For such $i, j$ and $a \leq Q_1$ we see that,
\beq
\R_{ij} \leq p^{-4/3},
\eeq
since $i \leq p / 200$ and $D_\S \leq K^{-\delta/2} p$.
This concludes the case $k=1$ using \eqref{eqn:USest}.  We assume that the estimate has been proven for $k-1$.  For $a \leq Q_k$, we have,
\beq
A_{ap}^{(k)} = \sum_{i, j} \US_{ai} \R_{ij} A_{jp}^{(k-1)} = \sum_i \sum_{j \leq Q_{k-1} } \US_{ai} \R_{ij} A_{jp}^{(k-1)}  +\sum_i \sum_{j>  Q_{k-1} } \US_{ai} \R_{ij} A_{jp}^{(k-1)}  =:G_1 +G_2.
\eeq
The induction assumption implies,
\beq
G_1 \leq \sum_{i, j} \US_{ai} \R_{ij} C p^{-4/3} \leq C p^{-4/3} ||\US \R||_{\ell^\infty \to \ell^\infty} \leq C p^{-4/3}.
\eeq
For the other term $G_2$, note that by the exponential decay we can restict the summation to $|i-a| \leq D_\S$.  Since $D_\S \leq K^{-\delta/2} C(Q_{k-1} - Q_k)$ we see that for such $i$ and $j  > Q_{k-1}$ that $\R_{ij} \leq C p^{-4/3}$.  Hence,
\beq
G_2 \leq \sum_{i, j } \US_{ai} Cp^{-4/3} A_{jp}^{(k-1)} \leq C p^{-4/3} ||\US||_{\ell^\infty \to \ell^\infty} ||A^{(k-1)} ||_{\ell^1 \to \ell^1} \leq Cp^{-4/3}.
\eeq
This concludes the proof. \qed

Therefore, we have the following estimate.
\bel \label{lem:fs}
Let $\delta_1 >0$.  
There is an event $\F_8$ of probability $\pp[\F_8] \geq 1 - K^{-c \delta_1}$ on which the following holds.  For all $p \geq K^{\delta_1} K^{5/6} \sqrt{1+t_1}$ and $a \leq K^{5/6}$,  we have
\beq
\UA_{ap} (s, t)  \leq C \frac{ (1+t_1)^{20}}{p^{4/3}}.
\eeq
for $t_0 \leq s \leq t \leq t_0 + t_1$.
\eel 
\proof Using the decomposition above we take $\eps_\ell = 1/10$ and $m=20$ so that $||B^{(m)}||_{\ell^\infty \to \ell^\infty} \leq C K^{-2}$.  The rest follows from the previous lemma. \qed

\subsection{Estimate for homogeneous solution}

We define $v(t)$ as the solution,
\beq
\frac{ \d }{ \d t } v(t) = \A (t) v (t) , \quad t_0 <t  < t_0 + t_1, \qquad v(t_0 ) = u (t_0),
\eeq
where $\A$ was defined above in \eqref{eqn:A-def}. 
We need to compare $v(t)$ to $u(t)$.
\bel \label{lem:uvest} Assume, 
\beq
\delta_c < 1/100. 
\eeq
There is an event $\F_9$ and a $c_1 >0$ with $\pp[\F_9] \geq 1- K^{-c_1}$ on which we have,
\beq
\sup_{t_0 \leq t \leq t_0 + t_1} \sup_{a \leq K^{1/2}} | u_a(t) - v_a (t) | \leq \frac{ (1+t_1)^{30}}{K^{1/50}}.
\eeq
for $N$ large enough.
\eel
\proof We have by the Duhamel formula, for any $\eps >0$,
\begin{align}
|u_a (t) - v_a (t) | & \leq \left| \sum_j \int_{t_0}^t \UA_{aj} (s, t) ( \xi_j (s) + q_j (s) \d s \right| \nonumber\\
&\leq \left| \sum_{j \leq K^{1-\eps}} \int_{t_0}^t \UA_{aj} (s, t) ( \xi_j (s) + q_j (s) \d s \right| + \left| \sum_{j > K^{1-\eps} } \sum_j \int_{t_0}^t \UA_{aj} (s, t) ( \xi_j (s) + q_j (s) \d s \right| \nonumber\\
& \leq \int_{t_0}^{t_0+ t_1} ||\1_{\{ j \leq K^{1-\eps} \}} \xi_j (s) ||_1 + || \1_{ \{ j \leq K^{1-\eps } \} } q_j (s) ||_1 \d s \nonumber \\
&+ \sum_{j > K^{1-\eps} } \int_{t_0}^{t_0+t_1}  |\UA_{aj} (s, t) (\xi_j (s) + q_j (s) ) | \d s =: A_1 + A_2.
\end{align}
Let us choose $\eps = 1/10$.  Let $\G' = \F_4 \cap \F_5$ be the intersection of the events of Lemma \ref{lem:xi} and \ref{lem:qi}, for some $\delta >0$ to be chosen.  Note that $A_1$ does not depend on the specific choice of $a \leq K^{1/2}$.  For $A_1$ we have,
\begin{align}
\ee[ \1_{ \G'} A_1 ] &\leq (1+t_1)^2 K^{1+2 \delta+\delta_c} \sum_{j \leq K^{1-\eps}} \frac{1}{ (K-j+1)^2} \nonumber\\
 \leq& C(1+t_1)^2 K^{1+2 \delta + \delta_c} \sum_{j \leq K^{1-\eps}} \frac{1}{K^2} \nonumber\\
\leq& C (1+t_1)^2 \frac{K^{2 \delta+\delta_c}}{K^{\eps}}.
\end{align}
Hence with probability at least $1-K^{-c}$ we have that $A_1 \leq (1+t_1)^2 K^{-1/20}$.  For $A_2$ we apply the estimate from Lemma \ref{lem:fs}.  Note that for $A_2$ we have $j \geq K^{9/10}$ and $a \leq K^{1/2}$ so Lemma \ref{lem:fs} is applicable.  Let $\G'' = \F_8 \cap \F_4 \cap \F_5 \cap \Gxi$ for sufficiently small $\xi >0$, where $\F_4$ and $\F_5$ are as before and $\F_8$ is the event of Lemma \ref{lem:fs} with $\delta_1 = 1/100$.  Then,
\beq
\1_{ \G''} A_2 \leq (1+t_1)^{29} \sum_{j > K^{1-\eps}} \frac{1}{j^{4/3}} \int_{t_0}^{t_1+t_0} \frac{ K^{1+\delta+\delta_c}}{ (K- j +1)^2} \frac{ B_{j, j+1} (s)}{j^{2/3}} \d s =: (1+t_1)^{30} A_3.
\eeq
Note that $A_3$ does not depend on the specific choice of $a \leq K^{1/2}$.  
We have,
\beq
\ee[ \1_{\Gxi} A_3] \leq (1+t_1)  K^{1+2 \delta+\delta_c} \sum_{j > K^{1-\eps}}^K \frac{1}{j^{4/3} (K-j+1)^2}.
\eeq
The sum we can estimate by,
\begin{align}
 \sum_{j > K^{1-\eps}}^K \frac{1}{j^{4/3} (K-j+1)^2} &= \sum_{K^{1-\eps} < j < K/2} \frac{1}{j^{4/3} (K-j+1)^2} + \sum_{ j > K/2} \frac{1}{j^{4/3} (K-j+1)^2} \nonumber\\
&\leq \frac{C}{K^2}  \sum_{K^{1-\eps} < j < K/2} j^{-4/3} + \frac{1}{K^{4/3}} \sum_{K/2 < j < K} \frac{1}{ (K-j+1)^2} \leq C K^{-4/3}.
\end{align}
Hence we get the claim by Markov's inequality and choosing, say $\delta =1/1000$. \qed

We require the following result which is  Proposition 10.4 of \cite{bourgade2014edge}.
\bep \label{prop:energy}
Let $\A$ be as above and consider the solution $\del_t w = \A w$.  Suppose that for some $b>0$, the coefficients of $\A$ satisfy,
\beq
B_{jk} \geq \frac{b}{ (j^{2/3}  - k^{2/3} )^2}
\eeq
and
\beq
W_j \geq \frac{b K^{1/3}}{ (K+1)^{2/3} - j^{2/3}}.
\eeq
Then for any $1 \leq p \leq  q \leq \infty$ and sufficiently small $\eta >0$, we have the estimate,
\beq
||w(t)||_q \leq C(q, p, \eta) \leq \left[ \left(K^{-2/3 \eta} t b \right)^{-(3/p-3/q) } \right]^{1-6\eta} ||v(0) ||_p.
\eeq
\eep
From this we obtain,
\bel \label{lem:vest}
Let $ \delta >0$.   Then there is a $C>0$ so that  on $\Gxi$ for $\xi >0$ sufficiently small we have, 
\beq
||v(t+t_0)||_\infty \leq C K^{\delta+\delta_c} t^{-3/20}.
\eeq
\eel
\proof For any $\delta_1 >0$ we have that,
\beq
||v(t_0)||_{10} \leq K^{\delta_1/10}
\eeq
on $\Gxi$ for sufficiently small $\xi$.  We can apply Proposition \ref{prop:energy} with $b = K^{-\delta_c}$, $p=10$, $q=\infty$ and $\eta \leq \min\{ \delta_1/10, 1/100\}$.   Hence, there is a constant $C>0$ so that,
\beq
||v(t+t_0)||_\infty \leq C K^{\delta_1} K^{ \delta_c} t^{-3/20}.
\eeq

\subsection{Proof of Theorem \ref{thm:te}}

We first prove the following.
\bep
Let $\delta >0$.  There is an event $\F'$ with probability at least $\pp[ F'] \geq 1- K^{-c \delta}$ on which we have for $t_1/2 \leq t \leq t_1$, the estimate
\beq \label{eqn:fin-est-1}
| x_i(t+t_0 ) - y_i (t+t_0) | \leq C \left(  \frac{(1+t_1)^{30}}{K^{1/50}} + \frac{K^{\delta+\delta_c}}{(t_1)^{3/20}} + (1+t_1) \frac{K^{\delta+\delta_c}}{K^{2/3}} + N^{-10} \right)
\eeq
for $i \leq K^{1/2}$.
\eep
\proof We write,
\begin{align}
x_i(t+t_0 ) - y_i (t+t_0) & = v_i (t+t_0) \label{eqn:tp-1} \\
&+ (u_i (t+t_0) - v_i (t+t_0) ) \label{eqn:tp-2} \\
&+ (\hatx_i (t+t_0) - \tilx_i (t+t_0) ) - ( \haty_i (t+t_0) - \tily_i (t+t_0) ) \label{eqn:tp-3}\\
&+ (x_i (t+t_0) - \hatx_i (t+t_0) ) - ( \haty_i (t+t_0 ) - y_i (t+t_0 ) )  \label{eqn:tp-4}
\end{align}
We now apply the previous results to estimate each of the above four terms.  From Lemma \ref{lem:vest} we have for \eqref{eqn:tp-1},
\beq
|v_i (t+t_0) | \leq \frac{K^{\delta_c+\delta}}{t^{3/20}}.
\eeq
on $\G_\xi$, for $\xi$ sufficiently small and all $0 \leq t \leq t_1$.  For \eqref{eqn:tp-2} we use Lemma \ref{lem:uvest}.  On the event $\F_9$ of that lemma we have, for $i \leq K^{1/2}$,
\beq
| u_i (t+t_0) - v_i (t+t_0) | \leq \frac{C (1+t_1)^{30}}{K^{1/50}}.
\eeq
for $0 \leq t \leq t_1$.  For \eqref{eqn:tp-3} on the event $\F_2$ of Lemma \ref{lem:cutoff-1} we have for $i \leq K^{1/2}$, 
\beq
| (\hatx_i (t+t_0) - \tilx_i (t+t_0) ) - ( \haty_i (t+t_0) - \tily_i (t+t_0) )| \leq (1+t_1) \frac{K^{\delta+\delta_c}}{K^{2/3}}.
\eeq
For \eqref{eqn:tp-4} we on the event $\F_1$ of Lemma \ref{lem:reg} that,
\beq
| (x_i (t+t_0) - \hatx_i (t+t_0) ) - ( \haty_i (t+t_0 ) - y_i (t+t_0 ) )  | \leq N^{-10}.
\eeq
This yields the claim, taking $\F' = \G_\xi \cap \F_1 \cap \F_2 \cap \F_9$. 
\qed

Now we complete the proof of the main theorem of this section.

\noindent{\bf Proof of Theorem \ref{thm:te}}.   We apply the previous proposition.  We take $K = N^{10^{-3}}$ and $t_1 = K^{10^{-6}}$.  For $\delta$ and $\delta_c$ sufficiently small, the RHS of \eqref{eqn:fin-est-1} is less than $CK^{-c}$ for some $c>0$.  The claim follows from taking $t_0 = -3 t_1/4$. \qed

\section{Properties of the limit} \label{sec:local}

\subsection{Locally Brownian properties: proof of Theorem \ref{thm:local-brown}}

\noindent{\bf Proof of \eqref{eqn:brown-inc}}. 
 We set $s= 0$ for notational simplicity.  The generalization is clear. 
By Theorem \ref{thm:conv-prob} there is an $1> \fa>0$ and $C_0 >0$ so that for all $N \geq 1$ we have $\pp [\F_N ] \geq 1 - C_0 N^{-\fa}$, where
\begin{align}
\F_N := \left\{ \sup_{ |t| \leq N^{\fa}, 1 \leq i \leq N^{\fa} } | \lambda_i (t) - \lamN_i (t) | \leq N^{-\fa} \right\}
\end{align}
For every $\delta >0$ there is a $C_\delta >0$ and an $\xi >0$ so that
\beq
\ee[ \1_{\F_\xi } |\lamN_{j+1} (t) - \lamN_{j} (t)|^{-1} ] \leq C_\delta N^{\delta} j^{1/3}, \qquad 1 \leq j \leq N^{1/2}.
\eeq
Let $\eps >0$ and $n >0$ as well as $T>0$.  Let $t_1 > t_0$ satisfy,
\beq
T \geq |t_1|, \quad T \geq |t_0|
\eeq
with $T \leq 1$.  
  Then,
\begin{align}
(\lambda_i(t_1 ) - \lambda_i (t_0) ) -(2/\beta)^{1/2} (B_i (t_1) - B_i (t_0 ) ) & = (\lambda^{(n)}_i(t_1 ) - \lambda^{(n)}_i (t_0) ) -(2/ \beta)^{1/2} (B_i (t_1) - B_i (t_0 ) ) \nonumber \\
&+ (\lambda^{(n)}_i(t_1 ) - \lambda^{(n)}_i (t_0) ) -(\lambda_i(t_1 ) - \lambda_i (t_0) ) 
\end{align}
Now the first term on the RHS can be written as 
\begin{align}
 (\lambda^{(n)}_i(t_1 ) - \lambda^{(n)}_i (t_0) ) - (2/\beta)^{1/2} (B_i (t_1) - B_i (t_0 ) )  &= \int_{t_0}^{t_1} \left( \sum_{j \neq i } \frac{1}{ \lambda^{(n)}_i (t) - \lambda^{(n)}_j (t) } - \frac{\lambda_i^{(n)} (t) }{2 n^{1/3}}\right) \d t \nonumber\\
&= \int_{t_0}^{t_1} \left( \sum_{j \neq i, |j-i| < n^{\delta_1} } \frac{1}{ \lambda^{(n)}_i (t) - \lambda^{(n)}_j (t) } \right) \d t  \nonumber\\
&+ \int_{t_0}^{t_1} \left( \sum_{|j-i | > n^{\delta_1} } \frac{1}{ \lambda^{(n)}_i (t) - \lambda^{(n)}_j (t) } - \frac{\lambda_i^{(n)} (t) }{2n^{1/2} }\right) \d t ,
\end{align}
where we choose
\beq
\delta_1 = \frac{ \fa }{10^3}.
\eeq
Let now, 
\beq
\E_{\xi, n} := \{  | \lambda_i^{(n)} (t) - n^{2/3} \gamma_i | < n^{\xi} (\hati)^{-1/3} : \forall i, \forall |t| \leq n \}.
\eeq
For every $\xi >0$ there is a $C_\xi$ so that,
\beq
\pp [ \E_{\xi, n} ] \geq 1 - C_\xi n^{-10}.
\eeq
On the event $\E_{\delta_2, n}$ with $\delta_2 = \delta_1 / 10$ we have for $n \geq C_2$ for some $C_2 >0$ that,
\beq
  \int_{-T}^T \left| \sum_{|j-i | > n^{\delta_1} } \frac{1}{ \lambda^{(n)}_i (t) - \lambda^{(n)}_j (t) } - \frac{\lambda_i^{(n)} (t) }{2n^{1/2} }\right| \d t \leq C_3 n^{10 \delta_2} (i)^{1/3}T
\eeq
On the other hand, by Markov's inequality, there is an event $\E_{1, i}$ and a $C_3 >0$ of probability at $1-n^{-\delta_1}- C_4 n^{-10}$ on which
\beq
\int_{-T}^{T} \left|  \sum_{j \neq i, |j-i| < n^{\delta_1} } \frac{1}{ \lambda^{(n)}_i (t) - \lambda^{(n)}_j (t) } \right| \d t \leq C_5 C_{\delta_1} n^{3 \delta_1}i^{1/3},
\eeq
for $i \leq n^{1/2}$, as long as $n \geq C_6$ for some $C_6 >0$. 
Therefore, on the event $\F_n \cap \E_{\delta_2, n} \cap \E_{1, i}$ which holds with probability at least
\beq
\pp [\F_n \cap \E_{\delta_2, n} \cap \E_{1, i} ] \geq 1- C_0 n^{-\fa} - n^{-\fa / 10^3} - C_4 n^{-10 \fa}
\eeq
we have for $n \geq C_2+C_6$ that,
\beq
\sup_{ -T \leq t_0  \leq t_1 \leq T }  \left| (\lambda_i(t_1 ) - \lambda_i (t_0) ) -(2 / \beta)^{1/2} (B_i (t_1) - B_i (t_0 ) ) \right| \leq n^{-\fa} + C_3 n^{\fa 10^{-3} }i^{1/3} T+ C_5 C_{\delta_1} n^{3 \fa/1000} i^{1/3} T,
\eeq
as long as $i \leq n^{\fa}$.
Choose now $T = \eps$.   Choose now $n = \lceil \eps^{-1/\fa} \rceil$.  We can assume $\eps \leq (C_2+C_6+100)^{-\fa}$ so that $n \geq C_2 + C_6+100$.  Then, for $i \leq \eps^{-1/100}$ we see that, (this implies $i \leq n^{\fa}$) 
\beq
\sup_{ -\eps \leq t_0  \leq t_1 \leq \eps }  \left| (\lambda_i(t_1 ) - \lambda_i (t_0) ) -(2 /\beta)^{1/2} (B_i (t_1) - B_i (t_0 ) ) \right|  \leq \eps + C'(\eps^{9/10} )
\eeq
for some $C' >0$ on an event with probability at least $1- C'' \eps^{10^{-4} }$.  This proves \eqref{eqn:brown-inc}. \qed

\noindent{\bf Proof of \eqref{eqn:brown-glob}}.  
We first define for notational simplicity,
\beq
f_j (s) := \lambda_j (s) - (2/\beta)^{1/2} B_j (s), \qquad \fn_j (s) := \lamn_j (s) -(2/\beta)^{1/2} B_j (s).
\eeq
For an integer $k$, let us define the event $\G_k$ as
\beq
\G_k := \left\{ \sup \left\{ |f_j (s) - f_j (t) | : |t|, |s| \leq 1, 2^{-k-1} \leq t-s \leq 2^{-k} \right\} >   2^{-kr}  \right\}.
\eeq
We estimate,
\beq
|f_j (s) - f_j (t) | \leq | \fn_j (s) - \fn_j (t) | + |f_j (s) - \fn_j (s) | + | \fn_j (t) - f (t) |.
\eeq
Now, 
\begin{align}
&\sup \left\{ | \fn_j (s) - \fn_j (t) | :  |t|, |s| \leq 1, 2^{-k-1} \leq t-s \leq 2^{-k}  \right\} \nonumber\\
\leq & \sup_{t \in \T_k } \int_{t}^{t+2^{-k+1}} \left| \sum_{ i \neq j} \frac{1}{ \lamn_j (s) - \lamn_i (s) } - \frac{ \lamn_j (s)}{2 n^{1/3}} \right| \d s 
\end{align}
where
\beq
\T_k := \{  t = j 2^{-k} : j \in \zz, |j| \leq 2^k \}.
\eeq
We have for any $p>0$ and $q>0$ and $i \leq n^{1/2}$ that
\beq
\ee \left( \int_{x}^y \frac{1}{ |\lamn_{i+1}(u) - \lamn_{i} (u) |} \d u \right)^{1+\beta-q} \leq C_{p, q} |x-y|^{\beta+1-q} i^{(1+\beta-q)/3} n^{p\fa}.
\eeq
Hence, arguing as above we see that, for $i\leq n^{1/2}$, and any $p>0, q>0$, $\eps >0$  and $D>0$ that,
\begin{align}
&\pp \left[ \int_{t}^{t+2^{-k+1}} \left| \sum_{ i \neq j} \frac{1}{ \lamn_j (s) - \lamn_i (s) } - \frac{ \lamn_j (s)}{2 n^{1/3}} \right| \d s > i^{1/3} n^{p \fa} 2^{-k} + \eps\right] \nonumber\\
 \leq & C' \left( n^{-\fa D} + \eps^{q-1-\beta} 2^{-k(\beta+1-q)} n^{p \fa} i^{(1+\beta-q)/3}\right)
\end{align}
for some $C'$ depending on $p, q$ and $D>0$ (but not on $\eps$).  We choose now $\eps = 2^{-k r}/4$.  Then, by the union bound,
\begin{align}
&\pp \left[ \sup_{t \in \T_k} \int_{t}^{t+2^{-k+1}} \left| \sum_{ i \neq j} \frac{1}{ \lamn_j (s) - \lamn_i (s) } - \frac{ \lamn_j (s)}{2 n^{1/3}} \right| \d s > i^{1/3} n^{p \fa} 2^{-k} + \eps\right]  \nonumber\\
\leq & C' ( n^{-\fa D} 2^k +  2^{-k(\beta+1-q) (1-r)+k } i^{(1+\beta-q)/3} n^{p\fa} )
\end{align}
Now, if 
\beq
r < 1 - \frac{1}{ 1+\beta}
\eeq
we can choose $q$ sufficiently small so that the exponent
\beq
\mfm = (\beta+1-q)(1-r) -1 > 0,
\eeq
is positive. 

Choose $n$ so that,
\beq
2^{k / \fa} \leq n \leq 2^{1+k/ \fa}.
\eeq
Choose $p$ so that,
\beq
p < \min \{ \mfa/10, (1-r)/10, \mfm / 10 \},
\eeq
and $D=100$.  
Hence, for $i$ such that
\beq
i \leq \min\{ 2^{k(1-r)/30} , 2^{k / (10 \fa) }, 2^{k \mfm / [10(1+\beta) ] } \}
\eeq
we have,
\begin{align}
&\pp\left[ \sup \left\{ | \fn_i (t) - \fn_i (s) | : |t|, |s| \leq 1, 2^{-k-1} \leq (t-s) \leq 2^{-k} \right\} > 2^{-kr}\left(\frac{1}{4} + C_1 2^{(r-1)k/2} \right)   \right]  \nonumber\\
\leq & C_2 \left( 2^{-k \mfm/10} + 2^{-k} \right)
\end{align}
Hence we see that there are exponents $\mfm_1 >0$ and $\mfm_2 >0$ so that, if $i \leq 2^{k \mfm_1}$ then,
\beq
\pp[ \G_k ] \leq C 2^{-k \mfm_2}.
\eeq
The claim follows from the union bound. \qed

\subsection{On the limiting SDE: proof of Theorem \ref{thm:sde}}

Note that Theorem \ref{thm:lr} implies the following estimate holds.
\beq \label{eqn:lr-est-2}
\pp^\beta[ \mu_{i+1}  - \mu_i \leq s i^{-1/3} ] \leq C \left( N^{\eps} s^{1+\beta-r} +N^{-D} \right)
\eeq
for any $r$ and $\eps >0$.
\bel
For any $K \geq \beta+1$ and $ r>0$ there is a constant $C_{K, r}$ so that  for every $i$ and $0<s<\frac{1}{4}$ satisfying,
\beq
s \leq \frac{1}{ i^{1/K}},
\eeq
so that the estimate,
\beq
\pp [ | \lambda_{i+1} (t) - \lambda_i (t) | \leq s i^{-1/3} ] \leq C_{K, r} s^{\beta+1-r},
\eeq
holds
\eel
\proof Let $0 < s < 1$.  Let $K>\beta+1$, $D>0$ and $\eps >0$ and $r>0$ be given.  Choose $n$ so that
\beq
\frac{1}{s^K} \leq n^{\fa} \leq \frac{2}{s^K}.
\eeq
By combining the estimate \eqref{eqn:lr-est-2} with Theorem \ref{thm:conv-prob} we have if $i \leq s^{-K}$ and $i^{1/3} \leq \frac{1}{2} s^{1-K}$, that,
\begin{align}
\pp\left[ |\lambda_{i+1} - \lambda_i |< s i^{-1/3} \right] \leq C \left( s^K + s^{-K\eps / \fa} s^{\beta+1-r} + s^{DK/\fa} \right).
\end{align}
Choosing $\eps $ so that $K\eps/ \fa = r$ and $D\geq \mfa $ large yields the claim.
\qed 

\bep \label{prop:sde-1}
Let $\frac{1}{10} > \om >0$ and let $K = N^{\om}$.  For $i \leq N^{\om /100}$ and for any $N \geq T>0, D>0$  there is a $C = C_{ \omega, D} >0$ so that,  
\beq
\pp\left[ \sup_{0 < t < T}\left| \int_0^t \left( \sum_{j > K } \frac{1}{ \lamN_i (s) - \lamN_j (s) } + N^{1/3} -  \left( \frac{ 16}{3 \pi^2} \right)^{1/3} ( \lfloor K \rfloor )^{1/3} \right) \d s \right| > T N^{-\omega/10} \right] \leq C N^{-D}
\eeq
\eep
\proof This is a straightforward consequence of rigidity.   Introduce the measure $\d \nu$ defined by,
\beq
\nu(E) \d E := N^{1/3 } \rhosc (N^{-2/3} E  ) \d E
\eeq
with quantiles $N^{2/3} \gamma_k$ where $\gamma_k$ are the $N$-quantiles of the semicircle distribution.  Let us introduce the shorthand,
\beq
\gamK := N^{2/3} \gamma_{ \lfloor K \rfloor}
\eeq
 Indeed on the event $\G_{\omega/100}$ we have for $N$ large enough,
\begin{align}
\left| \sum_{j > K } \frac{1}{ \lamN_i (s) - \lamN_j (s) } - \int_{\gamK}^\infty \frac{ \d \nu (E) }{-2N^{2/3} - E} \d E \right| &\leq C N^{\omega/100} \sum_{j > K } \frac{|i|^{2/3} + (\min \{ j^{1/3}, (N+1- j )^{1/3} \} )^{-1} }{ (i^{2/3} - j^{2/3} )^2} \nonumber\\
& \leq C N^{-\omega/20}
\end{align}
On the other hand, by the explicit formula of $\msc(z)$, 
\begin{align}
\int_{\gamK}^\infty \frac{ \d \nu (E) } {-2 N^{2/3}-E} = - N^{1/3}- \int_{-\infty}^{\gamK} \frac{ \d \nu (E) }{  -2 N^{2/3} - E}.
\end{align}
By the estimates,
\beq
\gamK = \left( \frac{2}{3} \lfloor K \rfloor \pi \right)^{2/3} - 2 N^{2/3} + \O \left( K^{5/3} N^{-1} \right)
\eeq
and for $-2N^{2/3} \leq E \leq \gamK$,
\beq
\nu (E) = N^{1/3}\rhosc(N^{-2/3} E ) = \frac{N^{1/3}}{ 2 \pi} \sqrt{ (2- E N^{-2/3} )(2+EN^{-2/3} ) } = \frac{1}{  \pi}  \sqrt{ 2 N^{2/3} +E }\left(1 + \O\left( |2+EN^{-2/3} | \right) \right)
\eeq
we see that,
\begin{align}
 \int_{-\infty}^{\gamK} \frac{ \d \nu (E) }{  -2 N^{2/3} - E} = - \left( \frac{ 16}{3 \pi^2} \right)^{1/3} ( \lfloor K \rfloor )^{1/3} + \O \left( K^2 N^{-2/3} \right).
\end{align}
This yields the claim. \qed

\bep \label{prop:sde-2}
Let $\fa$ be as in Theorem \ref{thm:conv-prob}.  Let $K = N^{\fa/100}$.  Let $i \leq K^{1/100}$.  There is a constant $C>0$ so that for any $N \geq T>0$ we have,
\begin{align}
\pp \left[ \int_{0}^T \left| \sum_{j \neq i, j \leq K } \frac{1}{ \lambda_i(s) - \lambda_j (s)} - \frac{1}{ \lamN_i (s) - \lamN_j (s) } \right| \d s > T N^{-\fa/4} \right]\leq C N^{-\fa/100}.
\end{align}
\eep
\proof Fix an $\eps = N^{-V}$ with $V$ to determined.  Let $\eps_{ij} = \eps$ for $i > j$ and $-\eps$ for $i < j$.  We estimate,
\begin{align}
\sum_{j \neq i, j \leq K } \left| \frac{1}{ \lambda_i(s) - \lambda_j (s)} - \frac{1}{ \lamN_i (s) - \lamN_j (s) } \right| & \leq \sum_{j \neq i, j \leq K } \left| \frac{1}{ \lambda_i(s) - \lambda_j (s) + \eps_{ij} } - \frac{1}{ \lambda_i (s) - \lambda_j (s) } \right| \nonumber \\
&+ \sum_{j \neq i, j \leq K } \left| \frac{1}{ \lamN_i(s) - \lamN_j (s) + \eps_{ij}} - \frac{1}{ \lamN_i (s) - \lamN_j (s) } \right| \nonumber \\
&+ \sum_{j \neq i, j \leq K } \left| \frac{1}{ \lambda_i(s) - \lambda_j (s) + \eps_{ij}} - \frac{1}{ \lamN_i (s) - \lamN_j (s) + \eps_{ij} } \right| \label{eqn:sde-bd}
\end{align}
By Markov's inequality for any $D>0$  there is a constant $C_{D}$ such that for any $\delta >0$ we have,
\beq
\pp \left[  \int_{0}^T \sum_{j \neq i, j \leq K } \left| \frac{1}{ \lamN_i(s) - \lamN_j (s) + \eps_{ij}} - \frac{1}{ \lamN_i (s) - \lamN_j (s) } \right|  > \delta \right] \leq C_{D}  N^{-D} +  C_1 T \eps^{1/2} K^3 \delta^{-1} 
\eeq
where $C_1$ does not depend on $V$ or $D$.  
The same estimate holds for first the term on the RHS of \eqref{eqn:sde-bd} (in fact, the $N^{-D}$ term is not present here).  For the final term, on the event of Theorem \ref{thm:conv-prob} we have,
\begin{align}
& \int_0^T \sum_{j \neq i, j \leq K } \left| \frac{1}{ \lambda_i(s) - \lambda_j (s) + \eps_{ij}} - \frac{1}{ \lamN_i (s) - \lamN_j (s) + \eps_{ij} } \right| \d s \nonumber \\
& \leq N^{-\fa} \left( \int_0^T \sum_{j \neq i , j \leq K } \frac{1}{ ( \lambda_i (s)- \lambda_j(s) + \eps_{ij} )^2} \right)^{1/2} \left( \int_0^T \sum_{j \neq i , j \leq K } \frac{1}{ ( \lamN_i (s)- \lamN_j(s) + \eps_{ij} )^2} \right)^{1/2}
\end{align}
By Markov's inequality,
\beq
\pp \left[ \int_0^T \sum_{j \neq i , j \leq K } \frac{1}{ ( \lambda_i (s)- \lambda_j(s) + \eps_{ij} )^2} > T N^{\fa/10} \right] \leq C_r N^{-\fa/20} K^2 \eps^{-r}
\eeq
for any $r>0$, 
and a similar estimate for the other term with the $\lamN_i(s)$.  It remains to take, e.g., $\delta = T N^{-\fa}$, $\eps = N^{-10\fa}$, $r=1/100000$ and $D = 10 \fa$.  \qed

\noindent{\bf Proof of Theorem \ref{thm:sde}}.  We take $N$  the smallest  integer satisfying 
\beq
K \leq N^{\fa/100} \leq 2K.
\eeq
Assume $T \leq K$. 
Then, for $ i  \leq K^{1/100}$ write,
\begin{align}
&\lambda_i (t) - \lambda_i (0)  - (2/\beta)^{1/2} (B_i (t) - B_i (0) ) - \sum_{j \neq i, j \leq K } \int_{0}^t \frac{1}{ \lambda_i (s) - \lambda_j (s) } - a K^{1/3} \d s \nonumber\\
=& ( \lambda_i (t) - \lambda_i (0)  - \lamN_i (t) + \lamN_i (0) ) \nonumber\\
+& \left( \int_{0}^t \sum_{j > K} \frac{1}{ \lamN_i (s) - \lamN_j (s) } +N^{1/3} - a K^{1/3} \right) \\
+& \left( \int_{0}^t \sum_{j \neq i, j \leq K} \frac{1}{ \lamN_i (s) - \lamN_j (s) } - \frac{1}{ \lambda_i (s) - \lambda_j (s) } \d s \right)
\end{align}
On the event of Theorem \ref{thm:conv-prob} the supremum of the first term on the RHS is $\O (N^{-\fa})$.  Propositions \ref{prop:sde-1} (with the choice of $\omega = \fa/100$ and $D= \fa$) and \ref{prop:sde-2} handle the last two lines. \qed

\appendix 

\section{Proof of level repulsion estimates} \label{a:lr}

We will follow very closely the proof of \cite{bourgade2014edge}.  For this we require some notation.  For positive integers $K$, we write the configuration space for $ \bx \in \rr^K$ and $\by \in \rr^{N-K}$ as,
\beq
( \lambda_1, \dots, \lambda_N ) = ( x_1, \dots, x_K, y_{K+1}, \dots y_N ) =: ( \bx, \by)
\eeq
For given $\by$, we denote by $\mu_\by ( \d \bx )$ the conditional distribution of $\bx$ given $\by$ which has the form,
\beq
\mu_\by ( \d \bx ) = \frac{1}{ Z_{\by}} \e^{ - \beta N H_\by (x) } \d x,
\eeq
where
\begin{align}
H_\by (x) &:= \frac{1}{2} \sum_{i \in I } V_\by ( x_i ) - \frac{1}{N} \sum_{i, j \in I, i< j} \log |x_j - x_i | \nonumber\\
V_\by (x) &:= V(N^{-2/3} x) - \frac{2}{N} \sum_{j \notin I} \log |x - y_j |
\end{align}
where $I = [[1, K]]$.  Here,  $V(x) = \frac{1}{2} x^2$ is the quadratic potential but the arguments are not very specific to this choice.  Additionaly, $Z_{\by}$ is the normalization constant for the conditional measure. We denote by $\R_K ( \xi)$ the set of ``good'' boundary conditions,
\beq
\R_K ( \xi ) := \{ \by : |y_k -N^{2/3} \gamma_k | \leq N^{\xi} \hat{k}^{-1/3}, k \notin I \}.
\eeq
From Theorem 3.2 of \cite{bourgade2014edge}, we see that for any $\delta >0$ and $K \geq N^{\delta}$ we have,
\beq
\pp^{\mu_\by} [ y_{K+1} - x_K \leq s K^{-1/3}  ] \leq C N^{C \xi } s^{1+\beta -r},
\eeq
for $\by \in \R_K ( \xi )$ and $\xi $ sufficiently smal.  Hence, we have the following.
\bel \label{lem:lr-1}
Let $\delta, \eps,  r >0$.  Then there is  $C$ so that for $i \geq N^{ \delta}$ and $\xi$ small enough we have,
\beq
\pp^\beta \left[ \Fxi \cap \{  | \lambda_{i+1}  - \lambda_i  | \leq s  i^{-1/3} \} \right] \leq C N^{\eps} s^{1+\beta-r}.
\eeq
\eel
Clearly, Theorem \ref{thm:lr} follows from Lemma \ref{lem:lr-1} as well as the following. 
\bel \label{lem:lr-2}  For all sufficiently small $\xi >0$ the following holds.
For $1 \leq K \leq N^{1/2}$ we have for $y \in \R_k ( \xi)$,
\beq
\pp^{\mu_\by } [ |\lambda_{K+1} - x_K | \leq s K^{-1/3} ] \leq C s^{\beta+1} (K^2 + K N^{C \xi } )^{\beta+1}.
\eeq
\eel
We first prove a weaker estimate.  We denote $\mu_{\by, 0}$ the measure,
\beq
\mu_{\by, 0} := (Z^*)^{-1} (y_{K+1} -x_K)^{-\beta} \mu_{\by}.
\eeq
That is, the term $(y_{K+1} - x_K )^\beta$ is dropped from the measure $\mu_\by$.
\bel \label{lem:lr-3}
Let $ y \in \R_{K} ( \xi)$.  Then,
\beq
\pp^{\mu_{\by}} [ y_{K+1} - x_K \leq s K^{-1} ] \leq Cs (K^2 + K N^{C \xi } )
\eeq
The same estimate holds for $\mu_{\by, 0}$.
\eel
\proof Set $y_+ := y_{K+1}$ and $y_- = y_+ - a$ where $a = N^{\xi} K^{-1/3}$.  We will decompose the configuration space according to how many particles lie in the interval $[y_-, y_+]$.  We denote this number by $n$.  For any $\varphi$ satisfying $0 \leq \varphi \leq c$ with $c\leq 1$ sufficiently small, we consider
\begin{align}
Z_\varphi &:= \sum_{n=0}^K \int_{ (-\infty, y_- )^{K-n} } \prod_{j=1}^{K-n} \d x_j \int_{(y_- , y_+ - a \varphi )^n} \prod_{ j = K-n+1}^K \d x_j \left[ \prod_{i, j \in I , i < j } (x_j - x_i )^\beta \right] \e^{ - N \frac{\beta}{2} \sum_{j \in I } V_\by (x_j ) } \nonumber \\
&= \sum_{n=0}^K (1- \varphi )^{n + \beta n (n-1)/2} \int_{ (-\infty, y_- )^{K-n} } \prod_{j=1}^{K-n} \d w_j \int_{(y_- , y_+ )^n} \prod_{ j = K-n+1}^K \d w_j  \left[ \prod_{i < j \leq K-n} (w_j - w_i)^\beta \right] \nonumber\\
&\times \left[ \prod_{K-n < i < j \leq K } (w_j - w_i)^\beta \right]\left[ \prod_{i \leq K-n} \prod_{j = K-n+1}^K (y_- + (1- \varphi )(w_j - y_- ) - w_i )^\beta \right] \nonumber \\
& \times \exp \left[ - \frac{N \beta}{2} \sum_{j \leq K-n} V_\by (w_j ) + \sum_{j > K-n} V_\by (y_- + (1-\varphi ) (w_j  - y_- ) ) \right],
\end{align}
where we have made the following change of variables,
\beq
w_j := x_j, \quad j \leq K - n, \qquad w_j := y_- + (1- \varphi )^{-1} (x_j - y_- ), \quad K-n+1 \leq j \leq K.
\eeq
We seek to lower bound $Z_\varphi$ in terms of $Z_{\varphi = 0}$.   We work with each term indexed by $n$ in the summand separately.  The interactions between $i \leq K-n$ and $j \geq K-n+1$ can be estimated by,
\beq
[y_- + (1- \varphi ) (w_j - y_- ) - w_i ]^\beta = [ (1-\varphi ) (w_j - w_i ) + \varphi (y_- - w_i ) ]^\beta \geq [ (1-\varphi ) (w_j - w_i ) ]^\beta
\eeq
for any $w_i \leq y_-$.  We now estimate the effect of scaling the potential $V_\by$.  Fix a parameter $M$ satisfying
\beq
M = N^{C \xi},
\eeq
for $C$ as in Lemma \ref{lem:Vder} below. 
We have,
\begin{align}
\exp \left[ - \frac{N \beta}{2} V_\by (y_-  + (1- \varphi ) (w_j -y_- ) ) \right] &= \exp \left[ - \frac{ N \beta}{2} V^*_\by (y_- + (1- \varphi ) (w_j - y_- )) \right] \nonumber\\
&\times \prod_{K+1 \leq k \leq K+M} (y_k - y_-  - (1- \varphi )(w_j - y_- ) )^\beta, \label{eqn:lr-proof-1}
\end{align}
where we defined
\beq
V^*_\by (x) := V(xN^{-2/3} ) - \frac{2}{N} \sum_{ k > K+M} \log |x- y_k |, \qquad x \in (y_-, y_+ ).
\eeq
Now, we have
\beq
\left| V^* (y_- + (1- \varphi ) (w_j - y_- ) )- V^* (w_j ) \right| \leq \max_{x \in [y_-, y_+] } | ( V^*_\by (x ) )' | a \varphi \leq CN^{-1}\varphi M.
\eeq
The derivative of $V^*_\by$ is estimated in Lemma \ref{lem:Vder} below.  Hence, we have
\beq
\exp \left[ - \frac{N \beta}{2} V_\by^* (y_- + (1-\varphi ) (w_j - y_- )  )\right] \geq \e^{ - C \varphi M} \exp\left[ - \frac{N \beta}{2} V_\by^* (w_j ) \right], \qquad j > K-n.
\eeq
Since for $k \geq K+1$, and $j \geq K-n+1$,
\beq
[y_k - y_- - (1- \varphi ) (w_j - y_-)]^\beta = [ (1- \varphi) (y_k - w_j ) + \varphi (y_k - y_- )]^\beta \geq  [(1- \varphi ) (y_k - w_j )]^\beta
\eeq
we conclude the estimate,
\begin{align}
Z_\varphi &\geq \sum_{n=0}^K (1- \varphi)^{n+\beta n (n-1)/2 + \beta n (K-n) + \beta n M } \e^{- C \varphi n N^{C  \xi } } \int_{ (-\infty, y_- )^{K-n} } \prod_{j=1}^{K-n} \d w_j \int_{ (y_- , y_+)^n } \prod_{j = K-n+1}^K \d w_j \nonumber\\
& \times \left[ \prod_{i < j \leq K } (w_j - w_i )^\beta \right] \e^{ - \frac{N \beta}{2} \sum_{j \leq K } V_\by (w_j ) }
\end{align}
Since $n \leq K$ we have the estimate,
\beq
(1- \varphi)^{n+\beta n (n-1)/2 + \beta n (K-n) + \beta n M } \e^{- C \varphi n N^{C  \xi } }  \geq (1- \varphi )^{C(K^2 + K N^{C \xi } )}.
\eeq
Therefore,
\beq
\frac{ Z_\varphi}{ Z_{\varphi  = 0 } } \geq (1 - \varphi )^{C K^2 + C K N^{C \xi } }.
\eeq
Choose $\varphi := s K^{-1/3} a^{-1} = s N^{- \xi}$.  Therefore,
\beq
\pp^{ \mu_\by } [ y_{K+1} - x_K \geq s K^{-1/3} ] \geq 1 - Cs (K^2 + K N^{C \xi  } ).
\eeq
The proof for $\mu_{\by, 0}$ is very similar.  The only difference is just that the $k = K+1$ factor is missing from \eqref{eqn:lr-proof-1} in the case that $j=K$.  This modification does not affect the proof. \qed

\noindent{\bf Proof of Lemma \ref{lem:lr-2}}.  Recall the definition of $\mu_{\by, 0}$. For brevity we set $X := y_{K+1} - x_K$.  We have,
\beq
\pp^{\mu_{\by} } [ X \leq s K^{-1/3} ] =\frac{  \ee^{\mu_{\by, 0} } [ \1_{ X \leq s K^{-1/3}} X^\beta ] }{ \ee^{ \mu_{\by, 0} } [ X^\beta ] }.
\eeq
Using Lemma \ref{lem:lr-3}, 
\beq
\ee^{\mu_{\by, 0} } [ \1_{ X \leq s K^{-1/3}} X^\beta ] \leq (s K^{-1/3} )^\beta C sK^{-1/3} (K^2 +  K N^{C \xi } ).
\eeq
By the same lemma,
\beq
\pp^{ \mu_{\by, 0} } [ X \geq c K^{-1/3} ( K^2 + K N^{C \xi } )^{-1} ] \geq 0.5,
\eeq
and so,
\beq
\ee^{ \mu_{\by, 0 } } X^\beta \geq c \left( \frac{ 1}{ K^{1/3} (K^2 + K N^{C \xi } ) } \right)^\beta.
\eeq
This concludes the lemma.
\qed

\noindent{\bf Proof of Theorem \ref{thm:lr}}.  The proof follows from the above lemma.  The estimates for $i \geq N^{\delta}$ some small $\delta>0$ are proven using Lemma \ref{lem:lr-1} and the estimate for the remaining particles follow from Lemma \ref{lem:lr-2}.  Taking $\delta$ sufficiently small depending on $\eps$ gives the result. \qed

The following is similar to Lemma C.2 of \cite{bourgade2014edge}. 
\bel \label{lem:Vder}
For sufficiently small $\xi >0$ and $M = N^{C \xi}$ with $C \geq 10$  we have for $K \leq N^{1/2}$  and $\by \in \R_K ( \xi)$ that,
\beq
\max_{x \in [y_-, y_+] } \left|  \frac{x}{2 N^{4/3}} + \frac{1}{N} \sum_{j > K + M } \frac{1}{  y_j - x } \right| \leq \log(N) \frac{(K+M)^{1/3}}{N}
\eeq
where,
\beq
y_+ = y_{K+1}, \qquad y_- = y_+ - N^{\xi} K^{-1/3}
\eeq
\eel
\proof We introduce the measure $\nu (E)$ by
\beq
\nu (E) \d E:= N^{1/3} \rhosc (N^{-2/3} E) \d E.
\eeq
It is straightforward to see that for $x \leq y_+$,
\begin{align}
\left|  \frac{1}{N} \sum_{j > K + M } \frac{1}{  y_j - x } -\frac{1}{N} \int_{\gamM} \frac{ \d \nu (E)}{ E-x} \right| \leq C \frac{K^{1/3}}{N}.
\end{align}
where,
\beq
\gamM := N^{2/3} \gamma_{\lceil K+M \rceil }.
\eeq
Since,
\beq
\frac{1}{N} \int \frac{ \d \nu (E) }{E -x } = \frac{1}{N^{2/3}} \msc (x N^{-2/3})
\eeq
we see that
\beq
\left| \frac{1}{N} \int \frac{ \d \nu (E) }{E -x }  + \frac{x}{2 N^{4/3} } \right| \leq C N^{-2/3} (2 + xN^{-2/3} )_-^{1/2} \leq C \frac{N^{\xi}}{N}.
\eeq
It remains to prove,
\begin{align} \label{eqn:Vder-1}
\left| \int_{-\infty}^{\gamM} \frac{ \d \nu (E) }{ E -x } \right| \leq  C \log(N) (K+M)^{1/3}
\end{align}
If $x \leq -2 N^{2/3}$ then,
\beq
\left| \int_{-\infty}^{\gamM} \frac{ \d \nu (E) }{ E -x } \right| \leq C  \int_{0}^{ \gamM + 2N^{2/3}} \frac{1}{ \sqrt{E}} \d E \leq C(K+M)^{1/3}.
\eeq
Fix  $\eta = N^{-10}$.  If $0 \leq x + 2 N^{2/3} < \eta /2 $, let,
\beq
b = - 2 N^{2/3} + 2 (x+2 N^{2/3} )
\eeq
\beq
\left| \int_{-\infty}^{\gamM} \frac{ \d \nu (E) }{ E -x } \right|  \leq \left| \int_{-\infty}^{b} \frac{ \d \nu (E) }{ E -x } \right|  + \left| \int_{b }^{\gamM} \frac{ \d \nu (E) }{ E -x } \right| 
\eeq
The second term is less than,
\beq
 \left| \int_{b }^{\gamM} \frac{ \d \nu (E) }{ E -x } \right|  \leq C  \int_{0}^{ \gamM + 2N^{2/3}} \frac{1}{ \sqrt{E}} \d E \leq C(K+M)^{1/3},
\eeq
where we used that the denominator is greater than $\frac{1}{2}(E + 2 N^{2/3} )$.   The first term is less than,
\beq \label{eqn:Vder-2}
\left| \int_{-\infty}^{b} \frac{ \d \nu (E) }{ E -x } \right|  = \left| \int_{-\infty}^{b} \frac{ (\nu (E) - \nu (x) ) \d E }{ E -x } \right| \leq C \sqrt{x+2N^{2/3}}
\eeq
where we used the estimate,
\beq
| \nu(E) - \nu (x) |\leq C \frac{|E-x|}{ \sqrt{x+2N^{2/3}}} 
\eeq
and that the integral is of length $2(x+2N^{2/3} )$.  If $x+2N^{-2/3} > \eta /2 $ we instead split the integral over the three regions, the central one centered at $x$ with length $\eta$, and the other two being from $-2N^{2/3}$ to the leftmost endpoint of the central region, $x-\eta/2$, and the third then from the rightmost endpoint of the central region, $x+\eta/2$ until $\gamM$.  The same calculation leading to \eqref{eqn:Vder-2} tells us that the central region contributes $C \eta / \sqrt{x}$.  For the other two regions their contribution is bounded by
\beq
C (\sup_{ E \leq \gamM } \nu(E)) \int_{\eta}^{N} \frac{1}{ E} \d E \leq C \log(N)(K+M)^{1/3}
\eeq
This completes the proof. \qed

\section{Stochastic continuity} \label{a:sc}

The following is proved in an identical manner to Theorem 11.5 in Chapter 11 of \cite{erdos2017dynamical}.
\bel \label{lem:asc-1}
There is a $C>0$ so that the following holds.  Suppose that for some particles $x_1 \leq x_2 \leq \cdots x_N$ the estimates
\beq \label{eqn:asc-1}
 \left| \frac{1}{N} \sum_{i=1}^N \frac{1}{ x_i - (E+\i \eta ) } - \msc (E+ \i \eta ) \right| \leq \frac{N^{\eps}}{N \eta}
\eeq
hold for all $|E| \leq 10$ and $10 \geq \eta \geq N^{\eps-1}$.  Assume also that $-2 - N^{\eps-2/3} \leq x_1$ and $2+ N^{\eps-2/3} \geq x_N$.  Then, for all $i$,
\beq
|x_i - \gamN_i | \leq  C \frac{N^{C\eps}}{ N^{2/3} \min\{ i^{1/3}, (N+1-i)^{1/3} \}}.
\eeq
\eel
We have also the following.
\bel \label{lem:asc-2}
There is a $C>0$ so that the following holds.  Suppose that
\beq \label{eqn:a-rig1}
|x_i - \gamN_i | \leq   \frac{N^{\eps}}{ N^{2/3} \min\{ i^{1/3}, (N+1-i)^{1/3} \}}.
\eeq
for some particles $x_1 < x_2 < \dots < x_N$.  Then, for all $|E| \leq 10$ and $10 \geq \eta \geq N^{C\eps-1}$ we have,
\beq
\left| \frac{1}{N} \sum_{i} \frac{1}{ x_i - (E + \i \eta ) } - \msc (E+\i \eta ) \right| \leq C \frac{N^{C \eps}}{N \eta}.
\eeq
\eel
\proof  Write $z = E+ \i \eta $.  We have, for $ \eta \geq C N^{\eps-1}$ that,
\begin{align}
\left| \frac{1}{N} \sum_{i} \frac{1}{ x_i - (E + \i \eta ) } - \msc (E+\i \eta ) \right|  \leq  CN^{\eps}  \sum_{i=1}^N \int_{ \gamN_{i-1}}^{\gamN_{i}} \frac{1}{ | x - z|^2} \frac{\rhosc (x) }{N^{2/3} \min\{ i^{1/3}, (N+1-i)^{1/3} \} } \d x
\end{align}
due to the fact that $|x_i -z | \geq c |x - z|$ for $x_i$ obeying \eqref{eqn:a-rig1}.   From the asymptotics,
\beq
\gamN_i +2 \asymp \frac{i^{2/3}}{N^{2/3}}, \qquad 2-  \gamN_i  \asymp  \frac{(N+1-i)^{2/3}}{N^{2/3}}
\eeq
and the explicit form of $\rhosc (x)$ we see that there is a $C>0$ so that for all $ x\in [ \gamN_{i-1} , \gamN_i]$ and all $i$,
\beq
\frac{\rhosc (x) }{N^{2/3} \min\{ i^{1/3}, (N+1-i)^{1/3} \} }  \leq \frac{C}{N}.
\eeq
Furthermore,
\beq
\int_{-2}^2 \frac{1}{ |x-z|^2} \d x \leq \frac{C}{\eta}.
\eeq
The claim follows. \qed

Consider now the function,
\beq
f_N (z, t) = \sum_{i=1}^N \frac{1}{ \lamN_i (t) - z}.
\eeq
for $\Im[z] > 0$.  By the Ito lemma,
\begin{align}
\d f_N (z, t) &=  \d M_N(z, t) + (\partial_z f_N )(z, t) f_N (z, t) \d t + ( \partial_z^2 f_N )(z, t) (\beta^{-1} - 2^{-1} ) \notag\\
+& \frac{1}{2N^{1/3}} \left( f_N (z, t) + z ( \partial_z f_N ) (z, t) \right).
\end{align}
where the martingale $\d M_N (z, t)$ is
\beq
\d M_N (z, t) := - \sum_i \frac{ \d B_i (t) }{ ( \lamN_i (t) - z)^2 }.
\eeq
By the BDG inequality for any time $t_0$ and $t$ we have for any $D>0$ and $\eps >0$ that
\beq
\pp \left[ \sup_{ 0 \leq s \leq t } \left|  \int_{t_0}^s \d M_(z, u ) \right| > N^{\eps} \sqrt{Nt} / \Im[z]^2 \right] \leq CN^{-D}.
\eeq
for some $C = C_{\eps, D}$.  From the fact that $f_N (z, t)$ is Lipschitz with constant less than $N / \eta $ on the domain $\Im[z] \geq \eta$ we deduce the following.
\bel \label{lem:asc-3}
For any $D>0$ the following holds.  There is a $C_D$ so that for any $t_0$ we have that the estimate,
\beq
\sup_{ |E| \leq 10, 10 \geq \eta \geq N^{-10}, t \leq N^{-100}} \left| f_N (E+\i \eta  , t_0 + t ) - f_N (E+\i \eta, t_0 ) \right| \leq N^{-40}.
\eeq
holds with probability at least $1-N^{-D}$.
\eel
With all of the above we can easily deduce Lemma \ref{lem:stoch-cont}. 

\noindent{\bf Proof of Lemma \ref{lem:stoch-cont}.}  Fix a grid of times $\T := \{ i/N^{100}: i \in \zz \} \cap [-N, N]$ and $D>0$.  Define $x_i(t) := \lamN_i (t)$, and
\beq
m_x (z, t) := \frac{1}{N} \sum_{i} \frac{1}{ x_i - z}
\eeq
and use the notation $z: = E + \i \eta$.  
By Theorem \ref{thm:rig} and Lemma \ref{lem:asc-2} we see that, with probability at least $1-C_D N^{-D}$ the estimate \eqref{eqn:asc-1} holds for any $\eps >0$ as well as,
\beq \label{eqn:asc-2}
\Im [ m_x(2+N^{\eps-2/3} + \i N^{-10}, t) ] + \Im [ m_x(-2-N^{\eps-2/3} + \i N^{-10}, t) ] \leq N^{-5}.
\eeq
For the above we just used $\Im [m_x (E+\i \eta ) ] \leq \eta (\min_i |x_i(t) - E | )^{-2}$, and that for $E = \pm (2 + N^{\eps-2/3})$ the minimum is greater than $N^{\eps/2-2/3}$ on $\F_{\eps/2}$.  By Lemma \ref{lem:asc-3} we see that we can extend this estimate to all $|t| \leq N$ on an event of probability at least $1-C_D N^{-D}$.   Assuming that $\F_{\eps/2}$ holds on every $t \in \T$ we now argue that $x_N (t) \leq 2 +N^{\eps-2/3}$ and $x_1 (t) \geq -2 - N^{\eps-2/3}$.  Indeed, since the sample paths are continuous and $\F_{\eps/2}$ holds for every $t \in \T$, if at any time $t$ we had that a particle hits either $E = \pm (2 + N^{\eps-2/3})$ then $\Im [ m_X ( \pm (2+ N^{\eps-2/3} + \i N^{-10} ), t) ] \geq N^9$ at this point which violates \eqref{eqn:asc-2}.

We conclude that the assumptions of Lemma \ref{lem:asc-1} hold for all $|t| \leq N$ and so we conclude the proof. \qed


\bibliography{mybib}{}

\begin{thebibliography}{10}

\bibitem{AH}
A.~Adhikari and J.~Huang.
\newblock Dyson {B}rownian motion for general beta and potential at the edge.
\newblock {\em Probab. Theory Related Fields}, pages 1--58, 2020.

\bibitem{ANvM}
M.~Adler, E.~Nordenstam, and P.~Van~Moerbeke.
\newblock The {D}yson {B}rownian minor process.
\newblock {\em Ann. Inst. Fourier}, 64(3):971--1009, 2014.

\bibitem{ABG}
R.~Allez, J.-P. Bouchaud, and A.~Guionnet.
\newblock Invariant beta ensembles and the {G}auss-{W}igner crossover.
\newblock {\em Phys. Rev. Lett.}, 109(9):094102, 2012.

\bibitem{AG}
R.~Allez and A.~Guionnet.
\newblock A diffusive matrix model for invariant {B}eta-ensembles.
\newblock {\em Electron. J. Probab.}, 18, 2013.

\bibitem{AGZ}
G.~W. Anderson, A.~Guionnet, and O.~Zeitouni.
\newblock {\em An introduction to random matrices}, volume 118.
\newblock Cambridge university press, 2010.

\bibitem{bourgade2018extreme}
P.~Bourgade.
\newblock Extreme gaps between eigenvalues of {W}igner matrices.
\newblock {\em preprint, arXiv:1812.10376}, 2018.

\bibitem{bourgade2012bulk}
P.~Bourgade, L.~Erd{\H{o}}s, and H.-T. Yau.
\newblock Bulk universality of general $\beta$-ensembles with non-convex
  potential.
\newblock {\em J. Math. Phys.}, 53(9):095221, 2012.

\bibitem{bourgade2014edge}
P.~Bourgade, L.~Erd{\H{o}}s, and H.-T. Yau.
\newblock Edge universality of $\beta$-ensembles.
\newblock {\em Comm. Math. Phys.}, 332(1):261--353, 2014.

\bibitem{bourgade2014universality}
P.~Bourgade, L.~Erd{\H{o}}s, and H.-T. Yau.
\newblock Universality of general $\beta $-ensembles.
\newblock {\em Duke Math. J.}, 163(6):1127--1190, 2014.

\bibitem{homogenization}
P.~Bourgade, L.~Erd{\H{o}}s, H.-T. Yau, and J.~Yin.
\newblock Fixed energy universality for generalized {W}igner matrices.
\newblock {\em Comm. Pure Appl. Math.}, 69(10):1815--1881, 2016.

\bibitem{CHH}
J.~Calvert, A.~Hammond, and M.~Hegde.
\newblock Brownian structure in the {KPZ} fixed point.
\newblock {\em preprint, arXiv:1912.00992}, 2019.

\bibitem{cusp2}
G.~Cipolloni, L.~Erd{\H{o}}s, T.~Kr{\"u}ger, and D.~Schr{\"o}der.
\newblock Cusp universality for random matrices, {II}: The real symmetric case.
\newblock {\em Pure and Applied Analysis}, 1(4):615--707, 2019.

\bibitem{C-KPZ}
I.~Corwin.
\newblock Kardar-{P}arisi-{Z}hang universality.
\newblock {\em Notices Amer. math. Soc.}, 63(3):230--239, 2016.

\bibitem{CH}
I.~Corwin and A.~Hammond.
\newblock Brownian {G}ibbs property for {A}iry line ensembles.
\newblock {\em Invent. Math.}, 195(2):441--508, 2014.

\bibitem{DMV}
D.~Dauvergne, M.~Nica, and B.~Vir{\'a}g.
\newblock Uniform convergence to the {A}iry line ensemble.
\newblock {\em preprint, arXiv:1907.10160}, 2019.

\bibitem{DV}
D.~Dauvergne and B.~Vir{\'a}g.
\newblock Basic properties of the {A}iry line ensemble.
\newblock {\em preprint, arXiv:1812.00311}, 2018.

\bibitem{dyson}
F.~J. Dyson.
\newblock A {B}rownian-motion model for the eigenvalues of a random matrix.
\newblock {\em J. Math. Phys.}, 3(6):1191--1198, 1962.

\bibitem{ES}
A.~Edelman and B.~D. Sutton.
\newblock From random matrices to stochastic operators.
\newblock {\em J. Stat. Phys.}, 127(6):1121--1165, 2007.

\bibitem{erdos2019matrix}
L.~Erd{\H{o}}s.
\newblock The matrix {D}yson equation and its applications for random matrices.
\newblock {\em preprint, arXiv:1903.10060}, 2019.

\bibitem{cusp1}
L.~Erd{\H{o}}s, T.~Kr{\"u}ger, and D.~Schr{\"o}der.
\newblock Cusp universality for random matrices {I}: Local law and the complex
  hermitian case.
\newblock {\em Comm. Math. Phys.}, pages 1--76, 2020.

\bibitem{localrelaxation}
L.~Erd{\H{o}}s, B.~Schlein, H.-T. Yau, and J.~Yin.
\newblock The local relaxation flow approach to universality of the local
  statistics for random matrices.
\newblock {\em Ann. Inst. Henri Poincar{\'e} Probab. Stat.}, 48(1):1--46, 2012.

\bibitem{localspectral}
L.~Erd{\H{o}}s and H.-T. Yau.
\newblock Universality of local spectral statistics of random matrices.
\newblock {\em Bull. Amer. Math. Soc.}, 49(3):377--414, 2012.

\bibitem{Gap}
L.~Erd{\H{o}}s and H.-T. Yau.
\newblock Gap universality of generalized {W}igner and $\beta$-ensembles.
\newblock {\em J. Eur. Math. Soc.}, 17(8):1927--2036, 2015.

\bibitem{erdos2017dynamical}
L.~Erdos and H.-T. Yau.
\newblock A dynamical approach to random matrix theory.
\newblock {\em Courant Lecture Notes in Mathematics}, 28, 2017.

\bibitem{FF}
P.~L. Ferrari and R.~Frings.
\newblock On the partial connection between random matrices and interacting
  particle systems.
\newblock {\em J. Stat. Phys.}, 141(4):613--637, 2010.

\bibitem{gorin2020universal}
V.~Gorin and V.~Kleptsyn.
\newblock Universal objects of the infinite beta random matrix theory.
\newblock {\em preprint, arXiv:2009.02006}, 2020.

\bibitem{GS1}
V.~Gorin and M.~Shkolnikov.
\newblock Multilevel {D}yson {B}rownian motions via {J}ack polynomials.
\newblock {\em Probab. Theory Related Fields}, 163(3-4):413--463, 2015.

\bibitem{GS2}
V.~Gorin and M.~Shkolnikov.
\newblock Interacting particle systems at the edge of multilevel {D}yson
  {B}rownian motions.
\newblock {\em Adv. Math.}, 304:90--130, 2017.

\bibitem{gorin2018stochastic}
V.~Gorin and M.~Shkolnikov.
\newblock Stochastic {A}iry semigroup through tridiagonal matrices.
\newblock {\em Ann. Probab.}, 46(4):2287--2344, 2018.

\bibitem{HP}
D.~Holcomb and E.~Paquette.
\newblock Tridiagonal models for {D}yson {B}rownian motion.
\newblock {\em preprint, arXiv:1707.02700}, 2017.

\bibitem{HL}
J.~Huang and B.~Landon.
\newblock Rigidity and a mesoscopic central limit theorem for {D}yson
  {B}rownian motion for general $\beta$ and potentials.
\newblock {\em Probab. Theory Related Fields}, 175(1-2):209--253, 2019.

\bibitem{KNT}
M.~Katori, T.~Nagao, and H.~Tanemura.
\newblock Infinite systems of non-colliding {B}rownian particles.
\newblock In {\em Stochastic analysis on large scale interacting systems},
  pages 283--306. Mathematical Society of Japan, 2004.

\bibitem{KO}
Y.~Kawamoto and H.~Osada.
\newblock Finite-particle approximations for interacting {B}rownian particles
  with logarithmic potentials.
\newblock {\em J. Math. Soc. Japan}, 70(3):921--952, 2018.

\bibitem{fixed}
B.~Landon, P.~Sosoe, and H.-T. Yau.
\newblock Fixed energy universality of {D}yson {B}rownian motion.
\newblock {\em Adv. Math.}, 346:1137--1332, 2019.

\bibitem{landonyau}
B.~Landon and H.-T. Yau.
\newblock Convergence of local statistics of {D}yson {B}rownian motion.
\newblock {\em Comm. Math. Phys.}, 355(3):949--1000, 2017.

\bibitem{edgeDBM}
B.~Landon and H.-T. Yau.
\newblock Edge statistics of {D}yson {B}rownian motion.
\newblock {\em Electron. J. Prob., to appear}, 2017.

\bibitem{O1}
H.~Osada.
\newblock Interacting {B}rownian motions in infinite dimensions with
  logarithmic interaction potentials.
\newblock {\em Ann. Probab.}, 41(1):1--49, 2013.

\bibitem{O2}
H.~Osada.
\newblock Interacting {B}rownian motions in infinite dimensions with
  logarithmic interaction potentials {II}: {A}iry random point field.
\newblock {\em Stochastic Process. Appl.}, 123(3):813--838, 2013.

\bibitem{OT}
H.~Osada and H.~Tanemura.
\newblock Infinite-dimensional stochastic differential equations arising from
  {A}iry random point fields.
\newblock {\em preprint, arXiv:1408.0632}, 2014.

\bibitem{PS}
M.~Pr{\"a}hofer and H.~Spohn.
\newblock Scale invariance of the {PNG} droplet and the {A}iry process.
\newblock {\em Journal of statistical physics}, 108(5-6):1071--1106, 2002.

\bibitem{Q}
J.~Quastel.
\newblock Introduction to {KPZ}.
\newblock {\em Current Developments in Mathematics}, 2011(1), 2011.

\bibitem{RRV}
J.~Ramirez, B.~Rider, and B.~Vir{\'a}g.
\newblock Beta ensembles, stochastic {A}iry spectrum, and a diffusion.
\newblock {\em J. Amer. Math. Soc.}, 24(4):919--944, 2011.

\bibitem{Sodin}
S.~Sodin.
\newblock A limit theorem at the spectral edge for corners of time-dependent
  {W}igner matrices.
\newblock {\em Int. Math. Res. Not. IMRN}, 2015(17):7575--7607, 2015.

\bibitem{tao2010random}
T.~Tao and V.~Vu.
\newblock Random matrices: Universality of local eigenvalue statistics up to
  the edge.
\newblock {\em Comm. Math. Phys.}, 298(2):549--572, 2010.

\bibitem{tao2011random}
T.~Tao and V.~Vu.
\newblock Random matrices: universality of local eigenvalue statistics.
\newblock {\em Acta Math.}, 206(1):127--204, 2011.

\bibitem{T}
L.-C. Tsai.
\newblock Infinite dimensional stochastic differential equations for
  {D}yson’s model.
\newblock {\em Probab. Theory Related Fields}, 166(3-4):801--850, 2016.

\end{thebibliography}


\begin{thebibliography}{9999}
\bibitem[EKYY]{EKYY} Erdos-Knowles-Yau-Yin.  Semicircle for general class






\end{thebibliography}
\bibliographystyle{abbrv}

\end{document}